\documentclass[a4paper]{article}

\usepackage{amsmath}
    \usepackage{amsxtra}
    \usepackage{amstext}
    \usepackage{amssymb}
    \usepackage{latexsym}

\input xy

\xyoption{all}

\xyoption{poly}

\newtheorem{prop}{Proposition}[section]
\newtheorem{thm}{Theorem}[section]

\newtheorem{lem}{Lemma}[section]

\begin{document}

\begin{center}

\large\textbf{SMOOTH PROJECTIVE  SYMMETRIC VARIETIES \\WITH PICARD
NUMBER ONE}

\

\

\large{ALESSANDRO RUZZI}

\

\small{Institut Fourier,\\
100 rue des Maths,\\
BP 74 38402, St Martin d'Heres, \\
FAX: 0476514478\\
ruzzi@mozart.ujf-grenoble.fr}

\

\end{center}

\begin{abstract}
\noindent We  classify the smooth projective symmetric $G$-varieties
with Picard number one (and $G$ semisimple). Moreover we  prove a
criterion for the smoothness of the simple (normal) symmetric
varieties whose closed orbit is complete. In particular we  prove
that, given a such variety $X$ which is not exceptional, then $X$ is
smooth if and only if an appropriate toric variety contained in $X$
is smooth.

\

\noindent \textit{keywords}: Symmetric varieties, Fano varieties.

\


\end{abstract}
\

\

 \

A Gorenstein normal algebraic variety $X$ over $\mathbb{C}$ is
called a Fano variety if the anticanonical divisor is ample.  The
Fano surfaces are classically called  Del Pezzo surfaces. The
importance of Fano varieties in the theory of higher dimensional
varieties is similar to the importance of Del Pezzo surfaces in the
theory of surfaces. Moreover  Mori's program predicts that every
uniruled variety is birational to a fiberspace whose general fiber
is a Fano variety (with terminal singularities).

Often it is useful to subdivide the Fano varieties in two kinds: the
Fano varieties with Picard number  one and the Fano varieties with
Picard number strictly greater than one. For example, there are many
results which give an explicit bound to some numerical invariants of
a Fano variety (depending on the Picard number and on the dimension
of the variety). Often there is an explicit expression for the Fano
varieties of Picard number one and another expression for the
remaining Fano varieties.

We are mainly interested in the  smooth projective spherical
varieties with Picard number one. The smooth toric (resp.
homogeneous) projective varieties with Picard number one are just
projective spaces (resp. G/P with G simple and P maximal parabolic).
We classify the smooth projective symmetric $G$-varieties whose
Picard group is isomorphic to $\mathbb{Z}$ (with the hypothesis $G$
semisimple). One can easily show that they are all Fano, because the
canonical bundle cannot be ample.

In \cite{Ru} there are some partial results regarding the
classification of smooth Fano  projective symmetric varieties with
Picard number greater than one and low rank. See also \cite{P} for a
classification of the  smooth projective horospherical varieties
with Picard number one.

Let $G$ be a connected reductive algebraic group over $\mathbb{C}$
and let $\theta$ be an involution of $G$. Let $H$ be a closed
subgroup of $G$ such that $G^{\theta}\subset H\subset
N_{G}(G^{\theta})$. We will say that    a normal equivariant
embedding of $G/H$ is a symmetric variety. Any symmetric variety can
be covered by open $G$-subvarieties which are simple, i.e. have one
closed $G$-orbit. In Section One we  recall some known facts about
symmetric varieties and we fix the notations. In Section Two we give
a criterion for the smoothness of the simple symmetric varieties
whose closed orbit is complete (see Theorem~\ref{ecc1} and
Theorem~\ref{smooth not ecc}). It is easily showed that every
complete symmetric variety is covered by simple open symmetric
$G$-subvarieties whose closed orbits are complete. We  generalize
the results of Timashev (see \cite{T1}) and Renner (see \cite{R}) on
the equivariant embeddings of reductive groups. The idea of the
proof is the following one. First we reduce to the affine  case,
though there are some technical problems when the variety is
exceptional (see \cite{dCP1} for a definition of exceptional
symmetric variety). Given an affine $X$, we define a toric variety
$\overline{T\cdot x_{0}}$ contained in $X$ and we prove that $X$ is
smooth if and only $\overline{T\cdot x_{0}}$ is smooth. In such case
we relate the combinatorial description of $\overline{T\cdot x_{0}}$
to the colored fan of $X$ (see \cite{Br3} for a definition of
colored fan).

In Section Three we use the previous criterion to classify the
smooth complete symmetric $G$-varieties with Picard number one and
$G$ semisimple (see Theorem~\ref{classif}). We do not classify
directly the projective smooth symmetric varieties with Picard
number  one because there is a combinatorial classification of the
complete symmetric varieties, but it is not so easy to say which
ones are projective. We prove that the smooth complete symmetric
varieties with Picard number  one (and $G$ semisimple) have at most
two closed orbits. Because of this fact we  easily prove that they
are all projective. Observe that all the simple complete symmetric
varieties are projective by a theorem of Sumihiro.

\section{Notations and
Background}\label{intro}

In this section we  introduce the necessary  notations. The reader
interested to the embedding theory of spherical varieties can see
\cite{Br3} or \cite{T2}. In \cite{V1} is explained such theory in
the particular case of the symmetric varieties. Let $G$ be a
connected reductive algebraic group over $\mathbb{C}$ and let
$\theta$ be an involution of $G$. Let $H$ be a closed subgroup of
$G$ such that $G^{\theta}\subset H\subset N_{G}(G^{\theta})$. We say
that   $G/H$ is a symmetric homogeneous variety. An equivariant
embedding of $G/H$ is the data of a $G$-variety $X$ together with an
equivariant open immersion $G/H\hookrightarrow X$. A normal
$G$-variety is called a spherical variety if it contains a dense
orbit under the action of an arbitrarily chosen Borel subgroup of
$G$. One can show that an equivariant embedding of $G/H$ is a
spherical variety if and only if it is normal (see \cite{dCP1},
Proposition 1.3). In this case we say that it is a  symmetric
variety.   We  say that a subtorus of $G$ is split if
$\theta(t)=t^{-1}$ for all its elements $t$. We say that a split
torus of $G$ of maximal dimension is a maximal split torus and that
a maximal torus containing a maximal split torus is maximally split.
One can prove that any maximally split torus is $\theta$ stable (see
\cite{T2}, Lemma 26.5). We fix arbitrarily a maximal split torus
$T^{1}$ and a maximally split torus  $T$ containing $T^{1}$. Let
$R_{G}$ be the root system of $G$ with respect to $T$ and let
$R_{G}^{0}$ be the subroot system of the roots fixed by $\theta$. We
define $R_{G}^{1}=R_{G}\, \backslash\, R_{G}^{0}$. We can choose a
Borel subgroup $B$ containing $T$ such that, if $\alpha$ is a
positive root in $R_{G}^{1}$, then $\theta(\alpha)$ is negative (see
\cite{dCP1}, Lemma 1.2). One can prove that $BH$ is dense in $G$
(see \cite{dCP1}, Proposition 1.3). Let $D(G/H)$ be the set of
$B$-stable prime divisors of $G/H$; its elements are called colors.
Since $BH/H$ is an affine open orbit, the colors are the irreducible
components of $(G/H)\ \backslash\ (BH/H)$. We say that a spherical
variety is simple if it contains one closed orbit. Let $X$ be a
simple symmetric variety with closed orbit $Y$. Let $D(X)$ be the
subset of $D(G/H)$ consisting of the colors whose closure in $X$
contains $Y$. We say that $D(X)$ is the set of colors of $X$.

To each prime divisor $D$ of $X$, we can associate the normalized
discrete valuation $v_{D}$ of $\mathbb{C}(G/H)$ whose ring is the
local ring $\mathcal{O}_{X,D}$. One can prove that $D$ is $G$-stable
if and only if $v_{D}$ is $G$-invariant, i.e. $v_{D}(s\cdot
f)=v_{D}(f)$ for each $s\in G$ and $f\in \mathbb{C}(G/H)$. Let $N$
be the  set of all $G$-invariant valuations of $\mathbb{C}(G/H)$
taking value in $\mathbb{Z}$ and let $N(X)$ be the set of the
valuations associated to the $G$-stable prime divisors of $X$.
Observe that each irreducible component  of $X\ \backslash\ (G/H)$
has codimension one because $G/H$ is affine. (To define $N(X)$ we do
not need that $X$ is simple). Let $S:=T/\,T\cap H\simeq T\cdot
x_{0}$, where $x_{0}=H/H$ denotes the base point of $G/H$. One can
show that the group $\mathbb{C}(G/H)^{(B)}/\mathbb{C}^{*}$ is
isomorphic to the character group $\chi(S)$ of $S$ (see \cite{V1},
\S2.3); in particular, it is a free abelian group. We define the
rank of $G/H$ as the rank of $\chi(S)$;  it is also equal to the
dimension of $S$. We can identify the dual group
$Hom_{\mathbb{Z}}(\mathbb{C}(G/H)^{(B)}/\,\mathbb{C}^{*},\mathbb{Z})$
with the group  $\chi_{*}(S)$ of one-parameter subgroups of $S$; so
we can identify $\chi_{*}(S)\otimes \mathbb{R}$ with
$Hom_{\mathbb{Z}}(\chi(S),\mathbb{R})$. The  restriction map to
$\mathbb{C}(G/H)^{(B)}/\mathbb{C}^{*}$ is injective over $N$ (see
\cite{Br3}, \S3.1, Corollaire 3), so we can identify $N$ with a
subset of $\chi_{*}(S)\otimes \mathbb{R}$. We say  that $N$ is the
valuation semigroup of $G/H$. Indeed, $N$ is the semigroup
constituted by  the vectors contained in the intersection of the
lattice $\chi_{*}(S)$ with an appropriate rational polyhedral convex
cone $\mathcal{C}N$, called  the valuation cone. For each color $D$,
we define $\rho(D)$ as the restriction of $v_{D}$ to $\chi(S)$. In
general, the map $\rho:D(G/H)\rightarrow\chi_{*}(S)\otimes
\mathbb{R}$ is not injective.

Let $C(X)$ be the cone in $\chi_{*}(S)\otimes \mathbb{R}$ generated
by $N(X)$ and $\rho(D(X))$. We can recover $N(X)$ from $C(X)$.
Indeed $N(X)$ consists of  the  primitive vectors of the
one-dimensional faces of $C(X)$ which intersect $N$ (we say that a
vector $v\in \chi_{*}(S)\backslash\{0\}$ is primitive if an equality
$v=av'$, where $v'$ is  a vector  in $\chi_{*}(S)$ and $a$ is a
positive integer, implies $a=1$). We denote by
$cone(v_{1},...,v_{r})$ the cone generated by the vectors
$v_{1},...,v_{r}$. Given a cone $C$ in $\chi_{*}(S)\otimes
\mathbb{R}$ and a subset $D$ of $D(G/H)$, we say that  $(C,D)$ is a
colored cone if: \vspace{-1 mm} \begin{enumerate}
\item[(i)] $C$  is generated by  $\rho(D)$ and by a
finite number of vectors in $N$;
\item[(ii)] the relative interior of $C$ intersects   $\mathcal{C}N$.
\end{enumerate}

\begin{prop}[see \cite{Br3}, \S 3.3, Th\'{e}or\`{e}me]
\label{col-cone} The map $X\rightarrow (C(X),D(X))$ is a bijection
from the set of simple  symmetric varieties to the set of colored
cones.
\end{prop}

Given a  symmetric variety $\widetilde{X}$ (not necessarily simple),
let $\{Y_{i}\}_{i\in I}$ be \vspace{0.1 mm} the set of $G$-orbits.
Observe that $\widetilde{X}$ contains a finite number of $G$-orbits,
thus $\widetilde{X}_{i}:=\{x\in \widetilde{X}\ |\  \overline{G\cdot
x}\supset Y_{i}\}$ is open in $\widetilde{X}$ and is a simple
symmetric variety whose closed orbit is $Y_{i}$. We define
$D(\widetilde{X})$ as the set $\bigcup_{i\in
I}D(\widetilde{X}_{i})$. The family
$\{(C(\widetilde{X}_{i}),D(\widetilde{X}_{i}))\}_{i\in I}$ is
\vspace{0.2 mm} called the colored fan of $\widetilde{X}$ and
determines completely $\widetilde{X}$. Moreover $\widetilde{X}$ is
complete if and only if $\mathcal{C}N\subset\bigcup_{i\in
I}C(\widetilde{X}_{i}) $ (see \cite{Br3}, \S 3.4, Th\'{e}orem\`{e}
2).

\begin{prop}[see \cite{Br3}, \S 3.4, Th\'{e}or\`{e}me 1]
\label{col-fan} A family $\{(C_{i},D_{i})\}_{i\in I}$ of colored
cones is the colored fan of a  symmetric variety if and only if
\begin{enumerate}
\item[(i)] given   a colored cone $(C',D')$ such that $C'$
is a face of $C_{i_{0}}$ and $D'$ is the set $\{D\in D_{i_{0}}|\,
\rho(D)\in C_{i_{0}}\}$ for an appropriate $i_{0}\in I$, then
$(C',D')$ belongs to $\{(C_{i},D_{i})\}_{i\in I}$;
\item[(ii)] the intersections of $\mathcal{C}N$ with the relative interiors of the
$C_{i}$ are pairwise disjoint.
\end{enumerate}
\end{prop}

We recall the description of the sets $N$ and $\rho(D(G/H))$. Before
that, we need to associate a root system  to $G/H$. The subgroup
$\chi(S)$ of $\chi(T^{1})$ has finite index, so we can identify
$\chi(T^{1} )\otimes\,\mathbb{R}$ with $\chi(S)\otimes\,\mathbb{R}$.
Because $T$ is $\theta$-stable, $\theta$ induces an involution of
$\chi(T)\otimes\,\mathbb{R}$ which we call again $\theta$. The
inclusion $T^{1}\subset T$ induces an isomorphism of
$\chi(T^{1})\otimes\,\mathbb{R}$ with the $(-1)$-eigenspace of
$\chi(T)\otimes\,\mathbb{R}$ under the action of $\theta$ (see
\cite{T2}, \S 26). Denoted by $W_{G}$ the Weyl group of $G$ (with
respect to $T$), fix arbitrarily a $W_{G}$-invariant inner product
$(\, \cdot ,\cdot)$ over $\chi(T)\otimes \,\mathbb{R}$ which
coincides with the product given by the Killing form over
$span_{\mathbb{R}}(R_{G})$. We denote again $(\, \cdot ,\cdot)$ the
restriction of this inner product to
$\chi(T^{1})\otimes\,\mathbb{R}$. Thus we can identify
$\chi(T^{1})\otimes\,\mathbb{R}$ with its dual
$\chi_{*}(T^{1})\otimes\,\mathbb{R}$.  The set $R_{G,\theta}:=
\{\beta-\theta(\beta)\ |\ \beta\in R_{G}\}\backslash\,\{0\}$ is a
root system in $\chi(S)\otimes\mathbb{R}$ (see \cite{V1}, \S 2.3
Lemme), which we call the restricted root system of $(G,\theta)$; we
call the non zero $\beta-\theta(\beta)$ the restricted roots. If $G$
is semisimple, then the rank  of $G/H$ is equal to the rank of
$R_{G,\theta}$, i.e. the dimension of
$span_{\mathbb{R}}(R_{G,\theta})$. Usually we denote by $\beta$
(respectively by $\alpha$)  a root of $R_{G}$ (respectively of
$R_{G,\theta}$); often we denote by $\varpi$ (respectively by
$\omega$)  a weight of $R_{G}$ (respectively of $R_{G,\theta}$). In
particular, we denote by $\varpi_{1},...,\varpi_{n}$ the fundamental
weight of $R_{G}$ (we have chosen the basis of $R_{G}$ associated to
$B$). Notice however that the weights of $R_{G,\theta}$ are weights
of $R_{G}$. The  involution $\iota:=-\varpi_{0}\cdot\theta$ of
$\chi(T)$ preserves the set of simple roots; moreover $\iota$
coincides with $-\theta$ modulo the lattice generated by $R_{G}^{0}$
(see \cite{T2}, p.169). Here $\varpi_{0}$ is the longest element of
the Weyl group of $R_{G}^{0}$.  We denote by
$\alpha_{1},...,\alpha_{s}$ the elements of the basis
$\{\beta-\theta(\beta)\, |\, \beta\in R_{G}$ simple$\}\backslash
\{0\}$   of $R_{G,\theta}$. Let $b_{i}$ be equal to $\frac{1}{2}$ if
$2\alpha_{i}$ belongs to $R_{G,\theta}$ and equal to one otherwise;
for each $i$ we define $\alpha_{i}^{\vee}$ as the coroot
$\frac{2b_{i}}{(\alpha_{i},\alpha_{i})}\alpha_{i}$. The set
$\{\alpha^{\vee}_{1} ,...,\alpha^{\vee}_{s}\}$ is a basis of the
dual root system $R^{\vee}_{G,\theta}$, namely the root system
composed by the coroots of the restricted roots. We call the
elements of $R^{\vee}_{G,\theta}$ the restricted coroots.  Let
$\omega_{1},...,\omega_{s}$  be the fundamental weights of $R
_{G,\theta}$  with respect to $\{\alpha_{1} ,...,\alpha_{s}\}$ and
let $\omega^{\vee}_{1},...,\omega^{\vee}_{s}$  be the fundamental
weights of $R^{\vee}_{G,\theta}$ with respect to
$\{\alpha^{\vee}_{1} ,...,\alpha^{\vee}_{s}\}$. Let $C^{+}$ be the
positive closed Weyl chamber of $\chi(S)\otimes\mathbb{R}$, i.e. the
cone generated by the fundamental weights of $R_{G,\theta}$ and by
the vectors orthogonal to $span_{\mathbb{R}}(R_{G,\theta})$. We call
$-C^{+}$ the negative Weyl chamber. Given a dominant weight
$\lambda$ of $G$, we denote by $V(\lambda)$ the irreducible
representation of highest weight $\lambda$.

If $G$ is semisimple and simply connected, then the lattice
generated by the restricted roots is isomorphic to the character
group of $T^{1}/(T^{1}\cap N_{G}(G^{\theta}))$. Moreover the weight
lattice of $R_{G,\theta}$ is $\{\varpi-\theta(\varpi)\,|\, \varpi\in
\chi(T)\}$ and can be identified with the character group of
$T^{1}/T^{1}\cap G^{\theta}$. We want to give another description of
this  lattice. We say that a dominant weight $\varpi\in \chi(T)$ is
a spherical weight if $V(\varpi)$ contains a non-zero vector fixed
by $G^{\theta}$ (in this definition we do not require that  $G$ is
semisimple). In this case, $V^{G^{\theta}}$ is one-dimensional and
$\theta(\varpi)=-\varpi$. Thus we can think $\varpi$ as a vector in
$\chi(S)\otimes\mathbb{R}$. One can show that  the lattice generated
by the spherical weights coincides with the weight lattice of
$R_{G,\theta}$. See \cite{CM}, Theorem 2.3 or \cite{T2}, Proposition
26.4 for an explicit description of the spherical weights. Such
description implies that the set of dominant weights of
$R_{G,\theta}$ is the set of spherical weights and that $C^{+}$ is
the intersection of $\chi(S)\otimes\mathbb{R}$ with the positive
closed Weyl chamber of the root system $R_{G}$. Sometimes we call
$\sum\mathbb{Z}^{+}\omega^{\vee}_{i}$ the set of spherical
coweights.

The set $N$ is equal to $-C^{+}\cap\chi_{*}(S)$; in particular it
consists of the lattice vectors of the rational polyhedral convex
cone $\mathcal{C}N=-C^{+}$. The set $\rho(D(G/H))$ is equal to
$\{\alpha^{\vee}_{1} ,...,\alpha^{\vee}_{s}\}$ (see \cite{V1}, \S
2.4, Proposition 1 and Proposition 2). Thus, given a simple
symmetric variety $X$, the set $\rho(D(X))$ consists of the
primitive generators of the one-dimensional faces of $C(X)$ not
contained in $-C^{+}$; indeed $(\alpha_{i},\alpha_{j}^{\vee})>0$ if
and only if $i=j$. For each $i$ the fibre
$\rho^{-1}(\alpha_{i}^{\vee})$ contains  at most two colors (see
\cite{V1}, \S2.4, Proposition 1). We say that a simple restricted
root $\alpha_{i_{0}}$ is exceptional if there are two distinct
simple roots $\beta_{i_{1}}$ and $\beta_{i_{2}}$ in $R_{G}$ such
that: i) $\beta_{i_{1}}-\theta(\beta_{i_{1}})=
\beta_{i_{2}}-\theta(\beta_{i_{2}}) = \alpha_{i_{0}}$; ii) either
$\theta(\beta_{i_{1}})\neq -\beta_{i_{2}}$ or
$\theta(\beta_{i_{1}})= -\beta_{i_{2}}$ and $(\beta_{i_{1}},
\beta_{i_{2}})\neq0$. In this case we say that also
$\alpha_{i_{0}}^{\vee}$, $\theta$ and all the equivariant embeddings
of $G/H$ are exceptional. If $\alpha_{i_{0}}$ is exceptional then
the simple roots $\beta_{i_{1}}$ and $\beta_{i_{2}}$ are uniquely
defined by the previous properties. If $G/H$ is exceptional, then
$\rho$ is not injective. Moreover, if $H=N_{G}(G^{\theta})$, then
$\rho^{-1}(\alpha_{i}^{\vee})$ contains two colors if and only if
$\alpha_{i}^{\vee}$ is exceptional. We say that $G/H$ contains a
Hermitian factor if the center of $[G,G]^{\theta}$ has positive
dimension.  If  $G/H$ does not contain a Hermitian factor then
$\rho$ is injective (see \cite{V1}, \S 2.4, Proposition 1). If
$\rho$ is injective, we denote by $D_{\alpha^{\vee}}$ the unique
color contained in $\rho^{-1}(\alpha^{\vee})$.

{\em Remark 1.}  One can suppose that $G$ is the product of a
connected, semisimple, simply connected group and of a central split
torus in the following way. Let $\pi:\widetilde{G}\rightarrow [G,G]$
be the universal cover of the derived group of $G$. We consider the
group $G'=\widetilde{G}\times (T^{1}\cap Z(G)^{0})$ and the map
$\pi':G'\rightarrow G$ defined by $\pi'(g,t)=\pi(g)\cdot t$. The
group $G'$ acts transitively on $G/H$ by $g\cdot x:=\pi'(g)\cdot x$
for each $x\in G/H$ and $g\in G'$. Thus $G/H$ is isomorphic to
$G'/(\pi')^{-1}(H)$. Moreover there is a unique involution $\theta'$
of $G'$ such that $\pi'\theta'=\theta\pi'$ (see \cite{S}). The group
$(\pi')^{-1}(H)$ contains and normalizes $(\pi')^{-1}(G^{\theta})$,
hence it normalizes also
$\pi^{-1}([G,G]^{\theta})^{0}=(\widetilde{G})^{\theta'}$. Thus
$(\pi')^{-1}(H)$ contains and normalizes $(G')^{\theta'}$, so
$G'/(\pi')^{-1}(H)$ is a $G'$-symmetric variety. We can also suppose
that the stabilizer of the base point $x_{0}$ in $Z(G)^{0}$ is
constituted by the elements of order two of $Z(G)^{0}$. Indeed, it
is sufficient to consider $\widetilde{G}\times T''$, where  $T''$ is
a torus such that $T''/\{t\in T''\, |\, t^{2}=id\}$ is isomorphic to
$(T^{1}\cap Z(G)^{0})\ /\ (T^{1}\cap Z(G)^{0}\cap H)$.

{\em In the following, unless explicitly stated, we assume that $G$
is as in Remark 1}.

We say that $\theta$ is decomposable if $G/G^{\theta}$ is the
product of (non trivial) symmetric homogeneous varieties, up to
quotient by a finite group. Otherwise we say that $\theta$ is
indecomposable. If $\theta$ is indecomposable and $G/H$ contains a
Hermitian factor, we  say that $G/H$ is Hermitian. If $\theta$ is
indecomposable, we have the following possibilities: i) $G$ is an
one-dimensional torus; ii) $G$ is  simple and iii) $G=\dot{G}\times
\dot{G}$ with $\dot{G}$ simple and $\theta$ defined by
$\theta(x,y)=(y,x)$. In the last case $G/G^{\theta}$ is isomorphic
to $\dot{G}$ and is not Hermitian. If $G$ is semisimple, then
$\theta$ is indecomposable if and only if $R_{G,\theta}$ is
irreducible. For a classification of the indecomposable involutions
see \cite{T2}, \S26. If $G/H$ is Hermitian, then $Z(G^{\theta})$ has
dimension one and $R_{G,\theta}$ has type $BC_{l}$, $C_{l}$, $B_{2}$
or $A_{1}$. Suppose $G$ semisimple (simply connected) and $\theta$
indecomposable, then $\rho^{-1}(\alpha_{i}^{\vee})$ contains two
colors  if and only if $G/H$ is Hermitian, $H=G^{\theta}$ and
$\alpha_{i}^{\vee}$ is short (see \cite{T2}, pages 177-178).

Let $X$ be a simple symmetric variety with closed orbit $Y$. One can
show that there is a unique affine $B$-stable  open set $X_{B}$ that
intersects $Y$ and is minimal for this property. Moreover the
complement $X\,\backslash\, X_{B}$ is the union of  the $B$-stable
prime divisors   not containing $Y$ (see \cite{Br3}, \S 2.2,
Proposition). The stabilizer $P$ of $X_{B}$ is a parabolic subgroup
containing $B$. There is a Levi subgroup $L'$ of $P$ such that the
$P$-variety $X_{B}$ is the product $R_{u}P\times Z$ of the unipotent
radical of $P$ and of an affine $L'$-spherical variety $Z$ (see
\cite{Br3}, \S 2.3, Th\'{e}or\`{e}me). We can choose $L'$ and $Z$ so
that $x_{0}$ belongs to $Z$. Let $L$ be the standard Levi subgroup
of $P$. If $X$ is not exceptional, we will see that $Z$ is a
symmetric variety and that $L'$ can be chosen equal to $L$. Notice
that $P$ is the stabilizer of $\bigcup_{D\in D(G/H)\setminus\,\,
D(X) }D$. \vspace{0.1 mm}Given an union $D_{i_{1}}\cup...\cup
D_{i_{s}}$ of colors, we denote its stabilizer by
$P(\{D_{i_{1}},..., D_{i_{s}}\})$ or by $P_{G}(\{D_{i_{1}},...,
D_{i_{s}}\})$. Given a root $\beta$, let $U_{\beta}$ be the
unipotent one-dimensional subgroup of $G$ corresponding to $\beta$.
Given $\mu\in \chi_{*}(T)\otimes\mathbb{Q}\equiv \chi
(T)\otimes\mathbb{Q}$, we denote by $P(\mu)$ the parabolic subgroup
of $G$ generated by $T$ and by the subgroups $U_{\beta}$
corresponding to the roots $\beta$ such that $(\beta,\mu)\ \geq 0$.
Given a subgroup $K$ of $G$, we use the notation $H_{K}$
(respectively $B_{K}$) for the intersection $K\cap H$ (respectively
$K\cap B$).

Now, we want to describe the groups $P=P(I)$, where $I\subset
D(G/H)$. Observe that there is a Levi subgroup $L'$ of $P$ such that
$P/H_{P}=R_{u}P\times L'/H_{L'}$. Indeed, there is a Levi subgroup
$L''$ of $P\times\mathbb{C}^{*}$ such that
$P\times\mathbb{C}^{*}/H_{P}\times\{\pm 1\}= R_{u}P\times (L''\cap
P)/H_{L''\cap P}\times \mathbb{C}^{*}/ \{\pm 1\}$,  because of
\cite{Br3}, \S 2.3, Th\'{e}or\`{e}me applied to  an equivariant
embedding $Z$ of $G \times \mathbb{C}^{*}/H\times\{\pm 1\} $ such
that $D(Z)\equiv D(G/H)\,\backslash \,I$ (we have a one-to-one
correspondence between $D(G/H)$ and $D(G \times
\mathbb{C}^{*}/H\times\{\pm 1\})$, which associates to each $D\in
D(G/H)$ its $\mathbb{C}^{*}$-span in $G \times
\mathbb{C}^{*}/H\times\{\pm 1\}$).

To describe $P(I)$, it sufficient to consider the case where $G$ is
semisimple, $H=G^{\theta}$, $ \theta $ is indecomposable and
$I=\{D\}$. Indeed $P(D_{i_{1}},..., D_{i_{s}})$ is equal to the
intersections of the $P(D_{i})$. Moreover if
$(G,\theta)=(G_{1},\theta)\times (G_{2},\theta)$ and $D$ is a color
of $G/H$, then there is either a $B_{G_{1}}$-stable divisor $D_{1}$
of $G_{1}/ G_{1}^{\theta}$ or a $B_{G_{2}}$-stable divisor $D_{2}$
of $G_{2}/G_{2}^{\theta}$ such that $D$ is the image of $D_{1}\,
\times G_{2}/G_{2}^{\theta}$, respectively of
$G_{1}/G_{1}^{\theta}\,\times D_{2}$, in $G/H$; moreover $P(D)$ is
$P_{G_{1}}(D_{1})\times G_{2}$, respectively $ G_{1}\times
P_{G_{2}}(D_{2})$. Notice that $D_{1}$, respectively $D_{2}$, may be
non prime.

We begin with some recalls  from \cite{V1}, \S2.3-\S3.4 (see also
\cite{T2}, pages 177-178). Suppose by simplicity that $G$ is
semisimple, write $D(G/G^{\theta})=\{D_{1},...,D_{m}\}$ and let
$\pi':G\rightarrow G/G^{\theta} $ be the canonical projection. For
any $i$, the ideal of $(\pi')^{-1}(D_{i})$ is principal because we
have supposed $G$ simply connected; let $f_{i}$ be a generator of
$(\pi')^{-1}(D_{i})$. Let $\mathcal{P}_{+}\subset \mathbb{C}[G]^{*}$
be the multiplicative monoid constituted by the commune eigenvectors
to $B$ (with respect to the left translation) and to $G^{\theta}$
(with respect to the right translation), normalized so that they
assume the value one on the unit element  of $G$. Let
$\mathcal{P}_{+}^{G^{\theta}}$ be the subset of $\mathcal{P}_{+}$
constituted by the $G^{\theta}$-invariant vectors. Up to
normalization,   the vectors $\{f_{i}\} $ forms a basis of
$\mathcal{P}_{+}$. We can index the colors so that $f_{i}$ belongs
to $\mathcal{P}_{+}^{G^{\theta}}$ if and only if $i\leq r$; in
particular the $G^{\theta}$-eigenvalue $\chi_{i}$ of $f_{i}$ is
trivial if and only if $i\leq r$. Here $r=m-2rank\
\chi(G^{\theta})$. We can also suppose that
$\chi_{r+j}=-\chi_{r+j+rank\ \chi(G^{\theta})}$ for each
$j=1,...,rank\ \chi(G^{\theta})$. The monoid
$\mathcal{P}_{+}^{G^{\theta}}$ is free with basis
$\{f_{1},...,f_{r},f_{r+1}f_{r+r\!an\!k\,\chi(G^{\theta})+1},...,f_{r+r\!an\!k\,\chi(G^{\theta})}f_{m}\}$.
The map, which associates to any $f\in \mathcal{P}_{+}$ its
$B$-weight $ \omega(f)$, gives an isomorphism between
$\mathcal{P}_{+}^{G^{\theta}}$ and the set of spherical weights. (We
remark that in \cite{V1} is used a lightly different map). Observe
that $P(D_{i})$ is equal to the stabilizer of the line
$\mathbb{C}f_{i}$.

Suppose now $\theta$ indecomposable, $G$ semisimple and $I=\{D\}$.
Let $L(\{D\})$ be the standard Levi subgroup of $P(\{D\})$. If
$G/G^{\theta}$ is not Hermitian, then $\mathcal{P}_{+}^{G^{\theta}}=
\mathcal{P}_{+}$ and $ \omega(f_{1}),...,\omega(f_{m})$  are the
fundamental weights of $R_{G,\theta}$; thus the stabilizer of
$D_{\alpha_{i}^{\vee}}$ is $P(\omega_{i} )$. In this case, the
colors of $G/H$ are in one-to-one correspondence with the colors of
$G/G^{\theta}$. More generally, if $D_{i}\in D(G/G^{\theta})$ is the
unique color in $\rho^{-1}(\rho(D_{i}))$, then $\omega(f_{i})$
belongs to $\mathcal{P}_{+}^{G^{\theta}}$ and $P(\{D_{i}\})$ is
$P(\omega(f_{i}) )$; on the other hand, $\pi^{-1}(\pi(D_{i}))$ is
$D_{i}$, so also $P(\{\pi(D_{i})\})$ is $P(\omega(f_{i}) )$ (here
$\pi$ is the canonical projection $G/G^{\theta}\rightarrow G/H$).
Moreover, $R_{L(\{D_{i}\})}$ is spanned by the simple roots $\{\beta
\}$ such that either $\beta =\theta(\beta )$ or  $(\beta
-\theta(\beta ))^{\vee}\neq\rho(D_{i})$.

Next, we consider the colors $D$ of $G/G^{\theta}$ such that
$\rho^{-1}(\rho(D))\neq \{D\}$. Supposing $\theta$ indecomposable,
there is an exceptional root $\alpha_{i_{0}}$ if and only if
$G/G^{\theta}$ is Hermitian and $R_{G,\theta}$ has type $BC_{l}$. In
this case, $H=G^{\theta}=N_{G}(G^{\theta})$ and $2\alpha_{i_{0}}$ is
a restricted root; in particular $\alpha_{i_{0}}^{\vee}$ is the
short simple restricted coroot and
$\rho^{-1}(\alpha_{i_{0}}^{\vee})$ contains two colors, namely
$D_{r+1}$ and $D_{r+2}$. Write
$\alpha_{i_{0}}=\beta_{i_{1}}-\theta(\beta_{i_{1}})=\beta_{i_{2}}-\theta(\beta_{i_{2}})$
as in the definition of exceptional root. The weight of
$f_{r+1}f_{r+2}$ is the $i_{0}$-th fundamental spherical weight,
namely $\omega_{i_{0}}=\varpi_{i_{1}}+\varpi_{i_{2}}$. Thus
$f_{r+1}$ and $f_{r+2}$ have weights respectively $\varpi_{i_{1}}$
and $\varpi_{i_{2}}$; the stabilizers of the corresponding colors
are respectively $P(\varpi_{i_{1}})$ and $P(\varpi_{i_{2}})$. We
have $\rho(D)=(\beta_{i_{1}}-\theta(\beta_{i_{1}}))^{\vee}=
(\beta_{i_{2}}-\theta(\beta_{i_{2}}))^{\vee} $; moreover the root
system of the standard Levi subgroup of $P(\varpi_{i_{1}})$
(respectively $P(\varpi_{i_{2}})$) is spanned by the simple roots
different from $\beta_{i_{1}}$ (respectively from $\beta_{i_{2}}$).
In particular, $R_{L(\{D\})}$ contains properly the root system
spanned by the simple roots $\{\beta \}$ such that either $\beta
=\theta(\beta )$ or  $(\beta -\theta(\beta ))^{\vee}\neq\rho(D)$.
Moreover, there is an automorphism of $G$, fixing $T$ and $B$, which
exchanges $P(\varpi_{i_{1}})$ with $P(\varpi_{i_{2}})$;   one can
show by the Iwasawa decomposition (see \cite{T2}, p.168), that this
automorphism fixes $G^{\theta}$, so it induces an automorphism of
$G/G^{\theta}$ exchanging $D_{r+1}$ with $D_{r+2}$. For the
exceptional case, see also \cite{dCS}, \S 4.

If $G/H$ is Hermitian  and  non exceptional, then the restricted
root lattice has index two in the lattice generated by the spherical
weights; the inverse image of $\alpha_{i}^{\vee}$ contains two
colors if and only if $H=G^{\theta}$ and $\alpha_{i}^{\vee}$ is
short (see \cite{T2}, page 177-178). Therefore  the weight of
$f_{r+1}f_{r+2}$ is $\omega_{1}=2\varpi_{1}$ if $R_{G,\theta}=A_{1}$
or $B_{2}$ and $\omega_{l}=2\varpi_{l}$ otherwise. The vectors
$f_{r+1}$ and $f_{r+2}$ have the same weight, namely $\varpi_{1}$ if
$R_{G,\theta}=A_{1}$ or $B_{2}$ and $\varpi_{l}$ otherwise. Thus the
colors corresponding to $f_{r+1}$ and $f_{r+2}$ have the same
stabilizer: respectively $P(\varpi_{1})$ if $R_{G,\theta}=A_{1}$ or
$B_{2}$ and $P(\varpi_{l})$ otherwise. Moreover these colors have
the same image in $G/N_{G}(G^{\theta})$. Hence the stabilizer of a
color $D$ of $G/G^{\theta}$ is equal to the stabilizer of the image
of $D$ in $G/N_{G}(G^{\theta})$. There is a unique simple root
$\beta $ such that  $(\beta -\theta(\beta
))^{\vee}=\rho(D_{r+1})=\rho(D_{r+2})$: it is $\beta_{1}$ if
$R_{G,\theta}=A_{1}$ or $B_{2}$ and $\beta_{l}$ otherwise. The root
system $R_{L(\{D_{r+1}\})}\equiv R_{L(\{D_{r+2}\})}$ is generated by
the simple roots different from $\beta$.

Now we list  the exceptional indecomposable involutions. We have
three possibilities: i) $(G,\theta)$ has type $AIII$ and
$G/G^{\theta}$ is isomorphic to $SL_{n+1}/S(L_{l}\times L_{n+1-l})$
(with $n\neq 2l-1$); ii) $(G,\theta)$ has type $DIII$ and $G/
G^{\theta} $ is isomorphic to $SO_{4l+2}/GL_{2l+1}$; iii)
$(G,\theta)$ has type $EIII$ and $G/G^{\theta}$ is isomorphic to
$E_{6}/D_{5}\times \mathbb{C}^{*}$. The possibilities for the
stabilizer of a color $D$ with $\rho(D)$ exceptional are,
respectively: i) $P(\varpi_{l})$ and $P(\varpi_{n+1-l})$; ii)
$P(\varpi_{2l})$ and $P(\varpi_{2l+1})$; iii) $P(\varpi_{1})$ and
$P(\varpi_{6})$.

\section{Smoothness of Symmetric Varieties}\label{smooth}
In this section we want to classify the smooth symmetric varieties.
Clearly we can reduce ourselves to studying the simple ones. Let $X$
be a simple symmetric variety with closed orbit $Y$. Recall that
$X_{B}\cong R_{u}P\times Z$, where $P=P(D(G/H)\setminus D(X))$ and
$Z$ is a $L'$-spherical affine variety for an appropriate Levi
subgroup $L'$ of $P$; we  denote by $L$  the standard Levi subgroup
of $P$.

\begin{lem}\label{L-thetastable}
If $\rho^{-1}( \alpha^{\vee} )$ is contained in $D(X)$ for each
exceptional coroot  $\alpha^{\vee}$ in $\rho(D(X))$, then the group
$L$ is $\theta$-stable and $L/H_{L}$ is a symmetric space; moreover
we can choose  $L'$ equal to $L$. In particular, $L$ is
$\theta$-stable if $X$ is not exceptional.
\end{lem}

{\em Proof of Lemma~\ref{L-thetastable}.} By the description of $P$
in the previous section, we known that $P=P(\omega^{\vee})$ with
$\omega^{\vee}$ a spherical coweight such that
$\iota(\omega^{\vee})=\omega^{\vee}$. Thus  $R_{L}$ is stabilized by
the involution $\iota$, because
$\beta_{1}-\theta(\beta_{1})=\beta_{2}-\theta(\beta_{2})$ if and
only if $\beta_{1}=\beta_{2}$ or $\beta_{1}=\iota(\beta_{2})$ for
any two simple roots $\beta_{1}$ and $\beta_{2}$ in $R_{G}^{1}$. In
addiction $R_{L}$ contains $R_{G}^{0}$, thus it is stabilized by
$\theta=-\varpi_{0}\iota$, where $\varpi_{0}$ is the longest element
of the Weyl group of $R_{G}^{0}$. Therefore $P\cap\theta(P)=L$,
$H_{L}=H_{P}$ and we can suppose $L=L'$. Notice that
$L^{\theta}\subset L\cap H\subset N_{L}(G^{\theta})\subset
N_{L}(L^{\theta})$. $\square$

Under the hypotheses of the previous lemma,  $R_{L}$ is generated by
$R_{G}^{0}$ and by the simple roots $\beta$ in $R_{G}^{1}$  such
that $\rho^{-1}((\beta-\theta(\beta))^{\vee})$ is contained in
$D(X)$ (see the description of $P$ in the previous section). In
particular, $R_{L,\theta}^{\vee}=\rho(D(X))$ if $G/H$ has not
Hermitian factors.

\begin{lem}\label{affine}
Let $I\subset D(G/H)$ and let $L'$ be a Levi subgroup of $P(I)$ such
that $P(I)/H_{P(I)}=R_{u}P\times L'/H_{L'}$. Then the homogeneous
variety $L'/H_{L'}$ is affine; in particular $H_{L'}$ is reductive.
\end{lem}

{\em Proof of Lemma~\ref{affine}.} Let $X'$ be an equivariant
embedding of $G \times \mathbb{C}^{*}/H\times\{\pm 1\} $ such that
$D(X)\equiv D(G/H)\,\backslash \,I$ (we have a one-to-one
correspondence between $D(G/H)$ and $D( G/H\times
\mathbb{C}^{*}/\{\pm 1\})$, which associates to each $D\in D(G/H)$
its $\mathbb{C}^{*}$-span in $G/H\times \mathbb{C}^{*}/\{\pm 1\}$).
Let $G'=G \times \mathbb{C}^{*}$, $H'=H\times\{\pm 1\}$ and
$B'=B\times\{\pm 1\}$. We can suppose $X'_{B'}\cong R_{u}P\times Z'$
where $Z'$ is a $( L'\times \mathbb{C}^{*})$-spherical variety.
First, suppose that $(G,\theta)$ is a product of indecomposable
exceptional involutions. Then $R_{u}P\times L'/H_{L'}\times
\mathbb{C}^{*}/\{\pm 1\}$ is the intersection of the two affine open
sets $X'_{B'}$ and $G'/H'$, thus it is affine (observe that the
$(L'\times \mathbb{C}^{*})$-stable colors do not intersect
$X'_{B'}$). Hence $L'/H_{L'}$ is affine and $H_{L'}$ is reductive.
When $\theta$ is non-exceptional, this lemma is immediately implied
by the   Lemma~\ref{L-thetastable}. The general case follows because
we can write $(G,\theta)=(G_{1},\theta)\times (G_{2},\theta)$ where:
i) $(G_{1},\theta)$ is a product of indecomposable exceptional
involutions; ii) $(G_{2},\theta)$ is non-exceptional; iii)
$H=G_{1}^{\theta}\times H_{G_{2}}$; iv) $L'=(L'\cap G_{1})\times
(L'\cap G_{2})$; v) $H_{L'}=H_{L'\cap G_{1}}\times H_{L'\cap
G_{2}}$. Indeed, these facts imply that $H_{L'}$ is reductive, so
$L'/H_{L'}$ is affine. $\square$

We are mostly interested in the smooth complete symmetric varieties,
so we can reduce ourselves to studying the simple affine varieties
with a $L'$-fixed point.

\begin{lem} \label{pto fisso}
If $Y$ is projective, then there is a unique point in $Z$ fixed by
$L'$.
\end{lem}

{\em Proof of Lemma~\ref{pto fisso}.} Notice that $Z\cap Y$ is one
point $x$ and  is $L'$-stable. Moreover $x$ is unique because every
affine spherical variety is simple. $\square$

The hypothesis of the previous lemma is equivalent to the condition
$dim\ C(X)$ $=rank\ G/H$. Now, we prove that $X$ is not smooth if
$L$ is not $\theta$-stable.

\begin{thm}\label{ecc1} Let $X$ be a simple smooth symmetric variety with
open orbit $G/H$ and suppose that the closed orbit $Y$ is
projective. Then $\rho(D(X))$ does not contain exceptional coroots.
In particular, $L$ is $\theta$-stable.
\end{thm}

\begin{prop}\label{ecc2} Under the hypotheses of the previous theorem,
$\rho^{-1}( \alpha^{\vee} )$ is contained in $D(X)$ for each
exceptional coroot  $\alpha^{\vee}$ in $\rho(D(X))$. In particular,
$L$ is $\theta$-stable.
\end{prop}

{\em Proof of Theorem~\ref{ecc1}.} The theorem is implied by the
previous proposition because, by the Theorem~\ref{smooth not ecc},
$\rho$ is injective over $D(X)$ for any smooth simple symmetric $X$
(which has  the property stated in the Proposition~\ref{ecc2}). We
remark that  we do not use the Theorem~\ref{ecc1} in the proof of
the Theorem~\ref{smooth not ecc}. $\square$

Now, we begin the proof of the Proposition~\ref{ecc2}. Observe that
$X$ is smooth if and only if $Z$ is smooth. Moreover $Z$ is smooth
if and only if it is a $L'$-representation. Indeed the affine
$L'$-variety $Z$ can be embedded in a $L'$-representation and one
can shift the fixed point in 0. Moreover, up to taking a
subrepresentation, we can identify $Z$ with this representation
because of Luna's fundamental lemma (see \cite{L}, II.2). To prove
the theorem, it is sufficient to show that there is not a
representation of dimension $dim\, L'/H_{L'}$ containing a vector
$v$ with stabilizer $H_{L'}$. Observe that we can work with a
quotient of $L'$ by a finite central subgroup $K$ contained in
$H_{L'}$. Indeed, if  $Z$ is a $L'$-representation with the previous
properties, $K$ acts trivially on $Z$, so $Z$ is a
$(L'/K)$-representation of dimension $(L'/K)/(H_{L'}/K)$ and
$Stab_{(L'/K)}(v)=(H_{L'}/K)$. Remark  that $L'/K$ is no longer a
subgroup of $G$, however we can restrict to study $Z$, ``forgetting"
the embedding of  $Z$ in $X$. We want reduce to study varieties
$L'/H_{L'}$ with the following property:

\begin{quote}[*]
\emph{There are an indecomposable exceptional involution
$(G',\theta')$ and a color $D$ such that $\rho(D)$ is exceptional,
$L'$ is a Levi subgroup of the stabilizer $P'=P_{G'}(D)$,
$H_{L'}=(L')^{\theta'}$ and $P'/(P')^{\theta'}\equiv R_{u}P'\times
L'/H_{L'}$.}
\end{quote}

\begin{lem}\label{simm+var*}
Up to quotient by a finite central group in $H_{L'}$, we can choose
$L'$   isomorphic to a product $\prod_{i=0}^{s} L_{i}$ so that
$L'/H_{L'}=\prod L_{i}/H_{L_{i}}$, where $L_{i}/H_{L_{i}}$ is a
symmetric variety if $s=0$ and satisfies [*] otherwise.
\end{lem}

{\em Proof of Lemma~\ref{simm+var*}.} First, we   reduce to the case
of an indecomposable involution. Write $(G,\theta)$ as $\prod
(G_{i},\theta)$, where $G_{0}\,\cap\, L$ is $\theta$-stable and the
$(G_{i},\theta)$ with $i>0$ are indecomposable exceptional
involutions such that $G_{i}\cap L$ is not $\theta$-stable.  Then
$G/H= G_{0}/H_{G_{0}}\times \prod_{i>0}G_{i}/G_{i}^{\theta}$. Indeed
$G_{i}^{\theta}=N_{G_{i}}(G_{i}^{\theta})$ for each $i>0$, thus $H$
is generated by $G^{\theta} $ and $H_{G_{0}}$. Thus it sufficient to
prove the lemma for the $G_{i}/G_{i}^{\theta}$ with $i>0$. Finally,
suppose  $\theta$ an indecomposable exceptional involution and write
$P=P(I)$ with $I\subset D(G/H)$. Let $D$ be a color in $I$ such that
$\rho(D)$ is exceptional. If $I$ contains at least two colors, let
$\widetilde{L}$ be the standard Levi subgroup of $P(I\setminus
\{D\})$,  let $\overline{L}$ be the quotient of $\widetilde{L}$ by
$Z([\widetilde{L},\widetilde{L}])$ and let $\overline{P}$ be the
quotient of $P(I)\cap \widetilde{L}$ by
$Z([\widetilde{L},\widetilde{L}])$. We have $P(I)=R_{u}P(I\setminus
\{D\})\times (P(I)\cap \widetilde{L})$ and $P(I) / P(I) ^{\theta}=
R_{u}P(I\setminus \{D\})\times (P(I)\cap\widetilde{L} )/(P(I)\cap
\widetilde{L}) ^{\theta} =R_{u}P(I\setminus \{D\})\times
\overline{P}/ \overline{P}^{\theta}$; moreover    $\widetilde{L}/
\widetilde{L}^{\theta }=\overline{L}/\overline{L}^{\theta}=
[\overline{L},\overline{L}]/[\overline{L},\overline{L}]^{\theta}\times
Z(\overline{L})/ Z(\overline{L})^{\theta}$. On the other hand, there
is a Levi subgroup $\overline{L}\, '$ of
$\overline{P}\cap[\overline{L},\overline{L}]$ such that
$\overline{P}/\overline{P}^{\theta}=R_{u}\overline{P}\times\overline{L}\,
'/(\overline{L}\, ')^{\theta} \times  Z(\overline{L}\,
')/Z(\overline{L}\, ')^{\theta}$. We can write
$([\overline{L},\overline{L}],\theta)=(L_{1},\theta)\times(L_{2},\theta)$,
with $(L_{1},\theta) $ non-exceptional and $(L_{2},\theta)$
indecomposable exceptional ($R_{[\overline{L},\overline{L}],\theta}$
is contained in $R_{G,\theta}$ and  $R_{G,\theta}$ contains one
exceptional root). We have $\overline{L}\, '/(\overline{L}\,
')^{\theta}= (\overline{L}\, '\cap L_{1})/(\overline{L}\, '\cap
L_{1})^{\theta}\times (\overline{L}\, '\cap L_{2})/(\overline{L}\,
'\cap L_{2})^{\theta}$, where $(\overline{L}\, '\cap
L_{1})/(\overline{L}\, '\cap L_{1})^{\theta}$ is a symmetric
variety. On the other hand $(\overline{P}\cap
L_{2})/(\overline{P}\cap L_{2})^{\theta}=R_{u}(\overline{P}\cap
L_{2})\times (\overline{L}\,'\cap L_{2})/(\overline{L}\,'\cap
L_{2})^{\theta}$ and $\overline{P}\cap L_{2}$ is \vspace{0.15 mm}
the stabilizer  of $D\cap (L_{2}/L_{2}^{\theta})$ in $L_{2}$. Thus
we can choose $L'$ as the product of $Z(\widetilde{L})^{0}$ with the
inverse image of $\overline{L}\, '$ in $\widetilde{L}$. $\square$

We know that given a semisimple group $G$, a symmetric variety $G/H$
and an irreducible $G$-representation $V$, then $V^{H}$ has
dimension at most one; moreover, if $dim\,V^{H}=1$, $V$ has
dimension at least $dim\, G/H+1$. We want to prove a similar result
for the varieties  satisfying $[*]$.  Such varieties have dimension
$dim\, G /( G )^{\theta} -dim R_{u}P$, because $P/H_{P}=R_{u}P\times
L'/H_{L'}$. By simplicity, we assume $(G,\theta)=(G',\theta')$ in
the Lemmas \ref{dim2}\,-\,\ref{colore excptional} ($(G',\theta')$ is
as in the definition of the property $[*]$). Explicitly $dim \,L' /(
L' )^{\theta}$ is equal, respectively, to $nl-l^{2}+l$ if
$G/G^{\theta}=SL_{n+1}/S(L_{l}\times L_{n+1-l})$, to $ 2l^{2}+l $ if
$G/ G^{\theta} =SO_{4l+2}/GL_{2l+1}$ and to 16 if
$G/G^{\theta}=E_{6}/D_{5}\times \mathbb{C}^{*}$. In particular
$L'/(L')^{\theta}$ has dimension at least two. We consider
separately the case where $dim\,L'/(L')^{\theta}=2$, because in such
case $L'/(L')^{\theta}$ is isomorphic to a symmetric variety.

\begin{lem}\label{dim2}
If $L'/H_{L'}$ has dimension two, then it is isomorphic to
$SL_{2}/SO_{2}$. Thus, a non-trivial  $L'$-representation $V$ with
$V^{H_{L'}}\neq0$ has dimension at least three. Moreover, there are
neither a smooth affine embedding of $SL_{2}/SO_{2}$ nor a smooth
affine embedding of $SL_{2}\times\mathbb{C}^{*}/SO_{2}\times
\{\pm1\}$.
\end{lem}

{\em Proof of Lemma~\ref{dim2}.} If $dim\,L'/H_{L'}=2$, then $G/H$
is isomorphic to $SL_{3}/S(L_{1}\times L_{2})$ and $P$ is either
$P(\varpi_{1})$ or $P(\varpi_{2})$; thus $L'$ is isomorphic to
$GL_{2} $. By the Lemma~\ref{affine}, $H_{L'}^{0}$ must be a maximal
torus, hence $L'/H_{L'}$ is a symmetric variety isomorphic either to
$SL_{2}/SO_{2}$ or to $SL_{2}/N_{SL_{2}}(SO_{2})$. Observe that a
spherical $SL_{2}$-representation has dimension at least three. In
this case the representation is $S^{2}\mathbb{C}^{2}$, which is an
embedding of
$SL_{2}\times\mathbb{C}^{*}/N_{SL_{2}}(SO_{2})\times\{\pm1\}$, where
$\mathbb{C}^{*}$ acts freely over $\mathbb{C}^{2}$. We have only to
prove that $L'/H_{L'}$ is isomorphic to $SL_{2}/SO_{2}$. Suppose by
contradiction that $H_{L'}$ is isomorphic to the almost direct
product $N_{SL_{2}}(SO_{2})\cdot\mathbb{C}^{*}id$ and let $\sigma$
be an involution of $GL_{2}$ such that $GL_{2}^{\sigma}$ is
$H_{L'}^{0}$. Let $T'$ be a maximally split  torus of $GL_{2}$ with
respect to $\sigma$ and let $S'=T'/ H_{T'}$, then
$\chi(T')/\chi(S')$ contains an element of order two. Indeed
$\chi(S')$ is generated by $2\beta$ and $\chi(T')$ contains $\beta$
(here $\beta$ is a root of $GL_{2}$). We can identify
$\mathbb{C}(L'/H_{L'})^{ (B_{L'})}/\mathbb{C}^{*} $ with $
\mathbb{C}(G/H )^{ (B) }/\mathbb{C}^{*}$, thus we can identify
$\chi(S )$ with $\chi(S')$ as subset of $\chi(B)$; moreover we can
identity both $\chi(T )$ and $\chi(T')$ with $\chi(B)$. Thus also
$\chi(T )/\chi(S )$ contains an element of order two. But $\chi(S)$
is spanned by $\varpi_{1}+\varpi_{2}$, so
$\chi(T)/\chi(S)=\mathbb{Z}(\varpi_{1}+\chi(S))\cong \mathbb{Z}$, a
contradiction. (Here $\varpi_{1}$ and $\varpi_{2}$ are the
fundamental weights of $SL_{3}$). $\square$.

\begin{lem}\label{L''} Suppose that $L'/H_{L'}$ satisfies $[*]$.
Then there is a subgroup $K$ of $L'$ isomorphic to $L \cap
\theta(L)$ in such way that $L^{\theta}\subset L \cap \theta(L)$ is
isomorphic to a subgroup $K'$ of $K^{\theta}$. Let $V$ be a
$L'$-representation which contains a vector $v$ fixed by
$K^{\theta}$. Then the stabilizer of $v$ in $[K,K]$ is equal either
to  $[K,K]$ or to $[K,K]':=[K,K]\cap K'$.
\end{lem}
Remark that $K/K'$ is a symmetric variety. We need this lemma only
when $dim\,L'/(L')^{\theta}>2$.

{\em Proof of Lemma~\ref{L''}.} Write $L'=gLg^{-1}$; we have
$L^{\theta}\subset P^{\theta}=(L')^{\theta}$, because
$H=G^{\theta}$.  When $G=E_{6}$, one can show, by a dimensional
count,  that $H_{L'}^{0}$ is isomorphic to $L \cap\, \theta(L)=
Spin_{8}\times SO_{2}\times\mathbb{C}^{*}$ (one have to use the
reductivity of $H_{L'}$ and the fact that $L^{\theta}\subset
H_{L'}\subset L'$). In the other cases we claim that      there is
an automorphism $\varphi$ of $L'$ which sends $gL^{\theta}g^{-1}$
onto $L^{\theta} $. Let $K$ be the image of $g(L \cap
\theta(L))g^{-1}$ by $\varphi$, then $K^{\theta}$ contains
$K':=L^{\theta}=\varphi(gL^{\theta}g^{-1})$. Now, we prove the \
claim. \ If $G/H$ \ $=$ \  $GL_{n+1}/GL_{l}\times GL_{n+1-l}$ (and
$P=P(\varpi_{l})$), let $\{e_{1},...,e_{n+1}\}$ be the canonical
basis of $\mathbb{C}^{n+1}$. Let $V_{1}=span\{e_{1},...,e_{l}\}$,
$W=span\{e_{l+1},...,e_{n+1-l}\}$ and
$V_{2}=span\{e_{n+2-l},...,e_{n+1}\}$. Let $\varphi$ be the
isomorphism from $V_{1}$ to $V_{2}$ such that
$\varphi(e_{i})=e_{n+2-i}$. We have $L= GL(V_{1} )\times GL(W\oplus
V_{2} )$, $L'=GL(V_{1} )\times GL(gW\oplus gV_{2} )$ and
$L^{\theta}= \Delta(V_{1},V_{2})\times GL(W)$, where
$\Delta(V_{1},V_{2})$ is the diagonal  in $GL(V_{1} )\times
GL(V_{2})$ (we have identified $V_{1}$ with $V_{2}$ by $\varphi$).
Similarly, $gL^{\theta}g^{-1}= \Delta(V_{1},gV_{2})\times GL(gW)$,
where $\Delta(V_{1},gV_{2})$ is the diagonal  in $GL(V_{1} )\times
GL(gV_{2})$ (we have identified $V_{1}$ with $gV_{2}$ by
$g\circ\varphi$). The previous facts follows because  $R_{u}P$ acts
trivially over $V_{1}$ and over $\mathbb{C}^{n+1}/V_{1}$. On the
other hand, $L^{\theta}$ stabilizes $g(V_{2}\oplus W)$ because it is
contained in $L'$. Let $\widetilde{g}$ be $g-Id$, then
$\widetilde{g}\mathbb{C}^{n+1}$ is contained in $V_{1}$ and
$gv=v+\widetilde{g}v$ for each $v$. Notice that  $L^{\theta}$
preserves the decomposition $V_{1}\oplus gW\oplus gV_{2}$. Indeed,
given any $h\in L^{\theta}$ and any $w\in W$,
$hgw=hw+h\widetilde{g}w\in (W\oplus V_{1})\cap g(W\oplus V_{2})=gW
$; we can proceed similarly for $gV_{2}$. Moreover, given any $h\in
L^{\theta}$ and $v\in V_{1}$, we have $hg\varphi
(v)-g\varphi(hv)=h\widetilde{g}\varphi(v)-\widetilde{g}\varphi(hv)\in
V_{1}\cap gV_{2}=0$, thus $L^{\theta}$ is contained in
$gL^{\theta}g^{-1}$. Therefore $L^{\theta}=gL^{\theta}g^{-1}$.
Finally, if $G/H=SO_{4l+2}/GL_{2l+1}$, then $L'$ is isomorphic to
$GL_{2l+1}$ and $ gL^{\theta}g^{-1}$ is isomorphic to $Sp_{2l}\times
\mathbb{C}^{*}\subset GL_{2l}\times \mathbb{C}^{*}$. It sufficient
to prove that $L^{\theta}$ is conjugated to $gL^{\theta}g^{-1}$ in
$L'\equiv GL_{2l+1}$. First, we conjugate by an element $k_{1}$ of
$L'\equiv GL_{2l+1}$, so that $k_{1} L^{\theta} k_{1}^{-1}$ is
contained in $GL_{2l}\times \mathbb{C}^{*}$. Next, we conjugated by
an element $k_{2}$ of $GL_{2l}$ so that $k_{2}k_{1} L^{\theta}
k_{1}^{-1}k_{2}^{-1}$ is $Sp_{2l}\times \mathbb{C}^{*}\equiv
gL^{\theta}g^{-1}$. The last affirmation of the lemma follows
because $[K,K]/[K,K]'$ is a symmetric variety. $\square$

\begin{lem}\label{colore excptional}
Let $(G,\theta)$ be an indecomposable exceptional involution, let
$D$ be a color such that $\rho(D)$ is exceptional, let $P$ be the
stabilizer of $D$ and let $L'$ be a Levi subgroup of $P$ such that
$P/H_{P} =R_{u}P\times L'/H_{L'}$. Suppose that $dim\, L'/H_{L'}>2$
and let $K$, respectively $K'$,  be as in the Lemma~\ref{L''}. Let
$V$ be a  non-trivial irreducible $L'$-representation such that
$V^{\,H_{L'}}\neq0$. Then:

\begin{enumerate}

\item[(i)] $dim V\geq dim L'/H_{L'}$ and $dim V\neq dim\, L'/H_{L'}+1$;
\item[(ii)] if $dim V= dim L'/H_{L'}$, then $V^{H_{L'}}$
has dimension one and $Stab_{L'}(v)$ is not reductive for any $v\in
V^{H_{L'}}\backslash\{0\}$;
\item[(iii)] if there is $v\in V$ with stabilizer $H_{L'}$, then  $dim V\geq dim L'/H_{L'}+2$.
\end{enumerate}

\end{lem}

{\em Proof of the Lemma~\ref{colore excptional}.} Observe that the
first two points of the lemma imply the third one. Indeed, if $V$ is
as in the last point, then $V$ cannot have dimension equal to $dim
L'/H_{L'}$  because $H_{L'}$ is reductive by the Lemma~\ref{affine}.
Notice that $K$ contains a maximal torus; moreover  there are two
parabolic subgroups $Q$ and $Q^{-}$ of $L'$, containing respectively
$B_{L'}$ and the opposite of $B_{L'}$, such that $Q\cap Q^{-}=K$ is
a Levi subgroup of both them.

For each $r$, let $\{e_{1},...,e_{r}\}$ be the canonical basis of
$\mathbb{C}^{r}$. It is sufficient to study the $L'$-representations
with dimension at most $dim\, L'/H_{L'}+1$. If
$G/H=GL_{n+1}/\,GL_{l}\times GL_{n+1-l}\, (=SL_{n+1}/S(L_{l}\times
L_{n+1-l}))$, then  $L'  $ is isomorphic to
$GL(\bigoplus_{i=1}^{n-l+1}\mathbb{C}e_{i})$ $\times$ $
GL(\bigoplus_{i=n-l+2}^{n+1}\mathbb{C}e_{i})$ in such way that $K$
corresponds to $GL(\bigoplus_{i=1}^{l}\mathbb{C}e_{i})\times
GL(\bigoplus_{i=l+1}^{n-l+1}\mathbb{C}e_{i}) \times
GL(\bigoplus_{i=n-l+2}^{n+1}\mathbb{C}e_{i})$ and $K'$ corresponds
to $\Delta\times GL(\bigoplus_{i=l+1}^{n-l+1}\mathbb{C}e_{i})$,
where $\Delta$ is the ``diagonal" in
$GL(\bigoplus_{i=1}^{l}\mathbb{C}e_{i})\times
GL(\bigoplus_{i=n-l+2}^{n+1}\mathbb{C}e_{i})$. Recall that $n>2$.
Suppose first that $[K,K] $ acts trivially over a non-zero vector of
$V^{H_{L'}}$; in this case
$SL(\bigoplus_{i=n-l+2}^{n+1}\mathbb{C}e_{i})$ acts trivially on
$V$. Using the fact that $dim\,Z(K')=2$, it is not difficult to
prove that $V^{[K,K]}$ is one-dimensional and that $Z(K')$ acts not
trivially on it, a contradiction. Hence, both the simple factors of
$[K,K]$ act not trivially on $V$ (otherwise there is $(g,g)\in
\Delta$ such that $(g,g)\cdot v\neq v$). The representation
$W:=\bigoplus_{i=1}^{n-l+1}\mathbb{C}e_{i} \otimes
(\bigoplus_{i=n-l+2}^{n+1}\mathbb{C}e_{i})^{*}$ has dimension equal
to $dim\,L'/H_{L'}$. If there is another $L'$-representation, not
isomorphic to $W$ as $[L',L']$-representation, over which both the
simple factors of $[K,K]$ acts non-trivially and which has dimension
at most  $dim\, W+1$,  then either $n=2l=4$ or $l=1$. In the first
case  such representation is isomorphic to
$W_{1}:=\bigwedge^{2}\bigoplus_{i=1}^{3}\mathbb{C}e_{i} \otimes
(\bigoplus_{i=4}^{5}\mathbb{C}e_{i})^{*}$ (as
$[L',L']$-representation); moreover it is isomorphic to $W$ as
$[K,K]$-representation. In the second case, such representation is
isomorphic to
$W_{2}:=\bigwedge^{n-1}\bigoplus_{i=1}^{n}\mathbb{C}e_{i} $ (as
$[L',L']$-representation) and we have $[K,K]=[K,K]'=
SL(\bigoplus_{i=2}^{n }\mathbb{C}e_{i})$, thus $W_{2}$ is
isomorphic, as $[K,K]$-representation, to the direct sum  of a
one-dimensional trivial representation with the dual of the standard
representation. The decomposition $W=(\bigoplus_{i=1}^{l
}\mathbb{C}e_{i}\otimes(\bigoplus_{i=n-l+2}^{n+1}\mathbb{C}e_{i})^{*})\oplus
(\bigoplus_{l+1}^{n-l+1}\mathbb{C}e_{i}\otimes
(\bigoplus_{i=n-l+2}^{n+1}\mathbb{C}e_{i})^{*})$ is
$[K,K]$-invariant; thus $W^{[K,K]'}$ is one-dimensional. One can
show that $R_{u}Q$ stabilizes any $v\in W^{[K,K]'}$; moreover
$dim\,Stab_{L'}(v)= dim\,R_{u}Q\cdot [K,K]'+1$. On the other hand,
$dim\,Z(K)=2$ and $Z(K)\cdot v$ is at plus, one-dimensional.
Therefore $T'':=Stab_{Z(K)}(v)$ is one-dimensional and the connected
component of $Stab_{L'}(v)$ is the semidirect product $R_{u}Q
\ltimes (T''[K,K]')^{0}$, so it is not reductive. Given any $U$
isomorphic, as $[L',L']$-representation, to $W$,  to $W_{1}$ or to
$W_{2}$, one can study in the same manner $U^{[K,K]'}$ and
$Stab_{L'}(v)$ for each $v\in U^{[K,K]'}\backslash \{0\}$.

If $G=SO_{4l+2}$, then $L'$ is isomorphic to $GL(\mathbb{C}^{2l+1})$
in such way that $K$ corresponds to
$GL(\bigoplus_{i=1}^{2l}\mathbb{C}e_{i})\times GL(
\mathbb{C}e_{2l+1})$ and $K'$ corresponds to
$Sp(\bigoplus_{i=1}^{2l}\mathbb{C}e_{i})\times GL(
\mathbb{C}e_{2l+1})$. The non-trivial irreducible
$L'$-representations with dimension lesser or equal to $2l^{2}+l+1$
are isomorphic, as $[L',L']$-representations, to
$\mathbb{C}^{2l+1}$, to $\bigwedge^{2}\mathbb{C}^{2l+1}$, to
$\bigwedge^{2l-1}\mathbb{C}^{2l+1}$ or to
$\bigwedge^{2l}\mathbb{C}^{2l+1}$. Let $V$ be a $L'$-representation
which contains a non-zero vector $w$ fixed by $K'$ and that is
isomorphic, as $[L',L']$-representation, to $\mathbb{C}^{2l+1}$. We
have $(\mathbb{C}^{2l+1})^{[K,K]}=
(\mathbb{C}^{2l+1})^{[K,K]'}=\mathbb{C}e_{2l+1}$; thus $V$ is
$\mathbb{C}^{2l+1}\otimes (\bigwedge^{2l+1}\mathbb{C}^{2l+1})^{*}$
and $w$  is a multiple of $e_{2l+1}\otimes v$, where $v$ is a
non-zero vector of $(\bigwedge^{2l+1}\mathbb{C}^{2l+1})^{*}$. The
stabilizer of $e_{2l+1}\otimes v$ in $L'$ is $R_{u}Q^{-}\ltimes
(SL(\bigoplus_{i=1}^{2l}\mathbb{C}e_{i})\times
GL(\mathbb{C}e_{2l+1}))$. If $H_{L'}$ fixes $e_{2l+1}\otimes v$,
then  $H_{L'}$ is isomorphic to a subgroup of
$SL_{2l}\times\mathbb{C}^{*}$ (because the restriction of the map
$Stab_{L'}(e_{2l+1})\rightarrow Stab_{L'}(e_{2l+1})/R_{u}Q^{-}$ to
$H_{L'}$ would be an isomorphism). However, one can show that there
is not a reductive subgroup of $SL_{2l}\times \mathbb{C}^{*}$,
containing $Sp_{2l}\times \mathbb{C}^{*}$, which has dimension equal
to $dim\, H_{L'}$, a contradiction. One can study in the same manner
the $L'$-representations isomorphic, as $[L',L']$-representations,
to $\bigwedge^{2l}\mathbb{C}^{2l+1}$. On the other hand,
$\bigwedge^{2}\mathbb{C}^{2l+1}$ is isomorphic, as
$[K,K]$-representation, to
$\bigwedge^{2}\bigoplus_{i=1}^{2l}\mathbb{C}e_{i}\oplus
\bigoplus_{i=1}^{2l}\mathbb{C}e_{i}$, hence
$(\bigwedge^{2}\mathbb{C}^{2l+1})^{[K,K]'}$ is one-dimensional. The
stabilizer  in $L'$ of any $v\in
(\bigwedge^{2}\mathbb{C}^{2l+1})^{[K,K]'}\setminus\{0\}$ contains
$R_{u}Q$; moreover $dim\,Stab_{L'}= dim\,R_{u}Q\cdot[K,K]'  +1$. On
the other hand, $T'':=Stab_{Z(K)}(v)$ is one-dimensional because
$Z(K)\cdot(v)$ is, at plus, one-dimensional. Therefore, the
connected component of $Stab_{L'}(v)$ is the semidirect product
$R_{u}Q \ltimes (T''[K,K]'))^{0}$, so it is not reductive. One can
study $\bigwedge^{2l-1}\mathbb{C}^{2l+1}$ in the same manner of
$\bigwedge^{2}\mathbb{C}^{2l+1}$ (notice that
$\bigwedge^{2l-1}\mathbb{C}^{2l+1}$ is the dual representation of
$\bigwedge^{2}\mathbb{C}^{2l+1}$).

If $L$ has type $E_{6}$, then, up to isogeny,  $L'$ is isomorphic to
$Spin_{10}\times \mathbb{C}^{*}$, in such way that $K$ corresponds
to $Spin_{8}\times SO_{2}\times\mathbb{C}^{*}$ and $K'$ corresponds
to $Spin_{7}\times\mathbb{C}^{*}$. We have already observed that, by
a dimensional count, one can show that $H_{L'}^{0}$ is isomorphic to
$K$. Observe that the hypothesis $V^{H_{L'}}\neq0$ implies that
$Z(L')$ acts trivially on $V$. The irreducible non-trivial
$Spin_{10}$-representations with dimension at most 17 are the
half-spin representations and the first fundamental representation
$\mathbb{C}^{10}$. One can easily show that none of such
representation contains a vector fixed by $H_{L'}$. $\square$

{\em Proof of Proposition~\ref{ecc2}}. We work with a quotient
$\widetilde{L}$ of $L'$ by a central subgroup in $H_{L'}$, so that
we can write $\widetilde{L}=\prod_{j=0}^{s} L_{j}$ as in the
Lemma~\ref{simm+var*}.   We denote by $H_{\widetilde{L}}$ the image
of $H_{L'}$ in $\widetilde{L}$. We want to reduce ourselves to the
case where $\widetilde{L}/\widetilde{L}^{\theta}$ is a product of a
variety satisfying $[*]$ and, eventually, a one-dimensional torus.
Suppose $Z$ smooth and write \vspace{0.3 mm} $Z=\bigoplus_{h=1}^{m}
V(\lambda_{h})$. Write $x_{0}=\sum_{h=1}^{m}v_{h}$ with $v_{h}\in
V(\lambda_{h})$, so each $\widetilde{L}\cdot v_{h}$ is dense in
$V(\lambda_{h})$. Moreover \vspace{0.3 mm}
$Z(\widetilde{L})^{0}\cdot v_{h}$ is contained in
$\mathbb{C}^{*}v_{h}$, so $dim\,\widetilde{L}\cdot v_{h}\leq$
$dim\,([\widetilde{L},\widetilde{L}]\cdot v_{h})+1$. Up to re-index,
we can suppose that $L_{1}$ does not act trivially on
$V(\lambda_{1})$. Write $\lambda_{i}=\sum\mu_{i,j}$, where
$\mu_{i,j}$ is a $L_{j}$-weight for each $i$ and $j$, so
$V(\lambda_{i})$ is isomorphic to $\bigotimes_{j}V(\mu_{i,j})$. Thus
$V(\mu_{i,j})$ contains a $H_{\widetilde{L}}$-fixed vector $v$ for
each $i$ and $j$; hence, by the Lemmas~\ref{dim2} and~\ref{colore
excptional}, $dim\, V(\mu_{i,j})\neq 2$ for each $j$. Let $I$ be the
set $\{j\,|\, dim V(\mu_{1,j})>1 \}$. If $I $ contains at least two
elements, then $dim\, V(\lambda_{1})=\prod_{j\in I} dim\,
V(\mu_{i,j})> \sum_{ j\in I}dim\, V(\mu_{i,j})+1 \geq \sum_{j\in I}
dim\, (L_{j}/H_{ L_{j}})+1 \geq dim\, \widetilde{L}v_{1}$, a
contradiction. Thus $V(\sum \mu_{1,j})= V( \mu_{1,1})\otimes
V(\sum_{j\neq1}\mu_{1,j})$, where $dim V(\sum_{j\neq1}\mu_{1,j})=1$.
We have two possibilities: either $V(\lambda_{1})$ is a smooth
affine embedding of $L_{1}\cdot v_{1}$ or there is a one-dimensional
torus $T_{1}\subset \prod_{i\neq1}L_{i}$, such that $V(\lambda_{1})$
is a smooth affine embedding of $(L_{1}T_{1})\cdot v_{1}$; in
particular the dimension of $V(\lambda_{1})$ is at most
$dim\,L_{1}/H_{L_{1}}+1$. This fact implies, by the
Lemma~\ref{colore excptional}, that $dim\,L_{1}/H_{L_{1}}$ is two.

Finally, suppose that $dim\,L_{i}/H_{L_{i}}=2$ for each $i>0$. By
the Lemma~\ref{dim2},  $dim
V(\lambda_{i})=dim\,[\widetilde{L},\widetilde{L}]v_{i}+1$ for each
$i$. Moreover we can suppose that, for each $i>0$, $L_{i}$ acts
non-trivially on $V(\lambda_{j})$ if and only if $i=j$; otherwise
$\widetilde{L}\cdot\sum_{1}^{m} v_{i}$ is not dense. We can identify
$L_{1}/H_{L_{1}}$ with $GL_{2}/(\mathbb{C}^{*})^{2}$ and
$V(\mu_{1,1})$ with $S^{2}\mathbb{C}^{2}$. Let $n\in
N_{SL_{2}}((\mathbb{C}^{*})^{2})\setminus (\mathbb{C}^{*})^{2}$,
then $n\cdot v_{1}=-v_{1}$.  The center of $\widetilde{L}$ acts with
one character over each $V(\lambda_{i})$; moreover
$Z(\widetilde{L})\cdot\sum_{1}^{m} v_{i}$ has dimension \vspace{0.3
mm} equal to $m$, because $\widetilde{L}\cdot \sum_{1}^{m} v_{i}$ is
dense in $Z$. Thus there is $t\in Z(\widetilde{L})$ such that $t\sum
v_{i}=-v_{1}+\sum_{i>1} v_{i}$, so  $nt\in H_{\widetilde{L}}$.
Therefore $H_{\widetilde{L}}$ is not generated by
$\widetilde{L}^{\theta}$ and $H_{ L_{0}}$, a contradiction.
$\square$

\textit{In the rest of this section, we always suppose that $L $ is
$\theta$-stable. } We choose, as basis of  the root system $R_{L}$
of  $(L,T)$, the basis associated to the Borel subgroup $B_{L} $.
Let $R_{L,\theta}$ be the restricted root system of $(L,\theta)$; it
is contained in $R_{G,\theta}$. We arbitrarily choose an order of
the connected components of the Dynkin diagram of $R_{L,\theta}$ and
we define $R_{L,\theta}^{j}$ as the subroot system of $R_{L,\theta}$
corresponding to the $j$-th connected component of the Dynkin
diagram. Let $\{\alpha_{1}^{j},...,\alpha_{l_{j}}^{j}\}$ be the set
of simple roots of $R_{L,\theta}^{j}$; we index them as in \cite{H}.
Let $l$ be the rank of $G/H$ and let $n$ be the semisimple rank of
$G$.  In particular we suppose that, if $R_{L,\theta}^{j}$ has type
$A_{l_{j}}$,  then $(\alpha_{h}^{j},\alpha_{k}^{j})\neq 0$ only if
$h$ belongs to $\{k-1,k,k+1\}$.  We prove the following theorem.

\begin{thm}\label{smooth not ecc} Let $X$ be a simple symmetric variety with
open orbit $G/H$. Suppose that the closed orbit $Y$ is projective
and that $\rho^{-1}( \alpha^{\vee} )$ is contained in $D(X)$ for
each exceptional coroot  $\alpha^{\vee}$ in $\rho(D(X))$. Then the
closure of $T\cdot x_{0}$ in $X$ is normal. Moreover $X$ is smooth
if and only if  this toric variety is smooth. This is equivalent to
the following conditions:
\begin{enumerate}
\item[(i)] The restricted root system $R_{L,\theta}$ has type $\prod_{i=1}^{p}A_{l_{i}}$ for appropriate
integers $p,l_{1},...,l_{p}$. Moreover the rank of $G/H$ is greater
or equal to $\sum_{i=1}^{p}(l_{i}+1)$.
\item[(ii)] The cone $C(X)$ is  generated by a basis of $Hom_{\mathbb{Z}}(\chi(S),\mathbb{Z})$.
\item[(iii)] We can index the dual basis
$\{\lambda_{1}^{1},...,\lambda_{l_{1}+1}^{1},\lambda_{1}^{2},...,\lambda_{l_{2}+1}^{2},...,
\lambda_{l_{q}+1}^{q}\}$ of $\chi(S)$ so that $(
\lambda_{h}^{j},(\alpha_{k}^{i})^{\vee})=\delta_{h,k}\delta_{i,j}$
and
$\frac{1}{l_{j}+1}((l_{j}+1)\lambda_{i}^{j}-i\lambda^{j}_{l_{j+1}})$
is the $i$-th fundamental weight of $R_{L,\theta}^{j}$.

\end{enumerate}
\end{thm}

Before  proving the theorem, we  make some remarks.

{\em Remarks 2.} i) We request that $X$ is normal. Indeed there is a
non normal equivariant completion $\widetilde{X}$ of a symmetric
space such that the closure of $T\cdot x_{0}$ in $\widetilde{X}$ is
smooth (see \cite{V1}, \S 7.2). See also \cite{T1}, $\S10$ for an
example of an affine non-normal embedding $\widetilde{Z}$ of a
symmetric space such that the closure of $T\cdot x_{0}$ is normal.
Moreover the open orbits, respectively of $\widetilde{X}$ and
$\widetilde{Z}$,  are not exceptional.

ii) Given  $X$ which verifies the conditions of the theorem, we can
almost describe the associated colored cone. The cone $C(X)$ is the
dual cone of $cone(\lambda_{1}^{1},..., \lambda_{l_{q}+1}^{q})$ and
$\chi(S)$ is the free abelian group with basis
$\{\lambda_{1}^{1},..., \lambda_{l_{q}+1}^{q}\}$. The set
$\rho(D(X))$ is constituted by the primitive vectors of the 1-faces
of $C(X)$ not contained in $-C^{+}$; moreover $D(X)$ contains
$\rho^{-1}(\{\alpha^{\vee}\}_{\alpha\in R_{L,\!\theta}\, simple})$.
Given any $\alpha^{\vee} $ in $\rho(D(X))\setminus
R_{L,\theta}^{\vee}$, then $\rho^{-1}(\alpha^{\vee})$ contains two
elements, while $\rho^{-1}(\alpha^{\vee})\cap D(X)$ contains one
element. Any of the two  colors in $\rho^{-1}(\alpha^{\vee})$ can be
a color of $X$. Remark that $\rho$ is injective over $D(X)$.

iii) Suppose that $L/H_{L}$ is a group $\dot{L}$, i.e. suppose that
$L=\dot{L}\times \dot{L}$ and that $\theta$ is defined by
$\theta(x,y)=(y,x)$. The first condition of the theorem means that
$\dot{L}$ is the quotient of a product $\Pi (SL_{l_{i}+1}\times
\mathbb{C}^{*})$ by a finite subgroup  (some integer $l_{i}$ can be
0). If all the conditions of the theorem are verified, then
$\dot{L}$ is isomorphic to the product $\Pi GL_{l_{i}+1}$ and $Z$ is
isomorphic to the product $\Pi M_{l_{i}+1}$, where $M_{l_{i}+1}$ is
the algebra of complex matrices of order $l_{i}+1$. Observe that the
product of matrices defines a product over $\Pi M_{l_{i}+1}$ which
extends the product over $\Pi GL_{l_{i}+1}$. This means that $\Pi
M_{l_{i}+1}$ is an algebraic monoid whose unit group is $\Pi
GL_{l_{i}+1}$.

iv) Timashev (see \cite{T1}) studied a particular class of
equivariant embeddings of reductive groups. These embeddings are not
necessarily normal. Here $G$ is a product $\dot{G}\times \dot{G}$
and the involution $\theta$ is defined by $\theta(x,y)=(y,x)$.
Timashev supposes that there is a faithful $\dot{G}$-linear
representation $V=\bigoplus V(\mu_{i})$ with the following
properties. The $\mu_{i}$ are all distinct and the map
$\dot{G}\rightarrow \mathbb{P}(\bigoplus
End_{\mathbb{C}}(V(\mu_{i})))$ is injective. He considers the
closure of $\dot{G}$ in $\mathbb{P}(\bigoplus
End_{\mathbb{C}}(V(\mu_{i})))$ and obtains a criterion for this
variety to be normal, respectively smooth. Observe that $G$ acts on
$\bigoplus End_{\mathbb{C}}(V(\mu_{i}))\subset End_{\mathbb{C}}(V)$
and the action coincides with the one over $V\otimes V^{*}$ through
the canonical isomorphism $End_{\mathbb{C}}(V)\cong V\otimes V^{*}$.
Also Renner (see \cite{R}) considers the case where $G/H$ is
isomorphic to a reductive group $\dot{G}$, but he studies the affine
(normal) equivariant embeddings of $\dot{G}$ with a $G$-fixed point.
He classifies the smooth ones in the case where the center of
$\dot{G}$ has dimension one and the derived group of $\dot{G}$ is
simple.

{\em Proof}. First we reduce ourselves to the affine case. Observe
that $X$ is smooth if and only if $Z$ is smooth. Let $U$
(respectively $U_{L}$) be the unipotent radical of $B$ (respectively
of $B_{L}$), then $\mathbb{C}[Z]^{U_{L}}$ is isomorphic to
$\mathbb{C}[X_{B}]^{U}$, so $\mathbb{C}(G/H)^{(B)}/\mathbb{C}^{*}$
is isomorphic to $\mathbb{C}(L/H_{L})^{(B_{L})}/\mathbb{C}^{*}$;  in
our case $T$ is contained in $L$ and we can identify both
$\mathbb{C}(G/H)^{(B)}/\mathbb{C}^{*}$ and
$\mathbb{C}(L/H_{L})^{(B_{L})}/\mathbb{C}^{*}$ with $\chi_{*}(S)$.
We can identify $\mathbb{C}[X_{B}]^{(B)}/\mathbb{C}^{*}$ with
$C(X)^{\vee}\cap\,\chi(S)$  and
$\mathbb{C}[Z]^{(B_{L})}/\mathbb{C}^{*}$ with
$C(Z)^{\vee}\cap\,\chi(S)$ (see \cite{Br3}, $\S3.2$); in particular
we can identify $C(X)\subset \chi_{*}(S)\otimes\mathbb{Q}$ with
$C(Z)\subset \chi_{*}(S)\otimes\mathbb{Q}$. Let $cl(T\cdot x_{0},X)$
(respectively $\overline{T\cdot x_{0}}$) be the closure of $T\cdot
x_{0}$ in $X$ (respectively in $Z$). To reduce to the affine case,
it is now sufficient to prove that $cl(T\cdot x_{0},X)$ is normal
(respectively smooth) if and only if   $\overline{T\cdot x_{0}}$ is
normal (respectively smooth). This fact holds because the
intersection of $cl(T\cdot x_{0},X)$ with $X_{B}$ is
$\overline{T\cdot x_{0}}$, so $\overline{T\cdot x_{0}}$ is open in
$cl(T\cdot x_{0},X)$. Moreover $cl(T\cdot x_{0},X)$ is covered by
the translates of $\overline{T\cdot x_{0}}$ with respect to the
action of $N_{G^{\theta}}(T^{1})$.

Observe that the combinatorial data $C(Z),D(Z),\chi(S)$ of $Z$ does
not change if we substitute  $L$ by a finite cover or by a quotient
by a normal subgroup contained in $H_{L}$. Thus we can substitute
$L$ with a group of simply connect type, as we have done with $G$ in
the Remark 1. The only problem is that the group obtained is no
longer a subgroup of $G$. On the other hand we can restrict to
consider $Z$, ``forgetting" $X$. \emph{In the following, we always
consider an affine embedding $Z$ of a symmetric homogeneous variety
$L/H_{L}$, where $L$ is as in the Remark 1.}

\begin{lem}\label{norm}
The closure of $T\cdot x_{0}$ in $Z$ is normal.
\end{lem}

{\em Proof of Lemma~\ref{norm}.} We have a commutative diagram

\[ \xymatrix{ & \mathbb{C}[Z]^{(T)}/\mathbb{C}^{*}
\ar@{->>}[dr]& & &\\ \ \ \ \ \ \ \ \ \ \
 \mathbb{C}[Z]^{(B_{L})}/\mathbb{C}^{*}\ar@{^{(}->}[ru] \ar@{^{(}->}[rr]&
 &\mathbb{C}[\overline{T\cdot x_{0}}]^{(T)}/\mathbb{C}^{*}.
 &   (2.1)  }   \]

Let $W_{L,\theta}$ be the Weyl group of $R_{L,\theta}$; it is
isomorphic to $N_{L^{\theta}}(T^{1})/C_{L^{\theta}}(T^{1})$
\vspace{0.3 mm} (see \cite{T2}, Proposition 26.2). Observe that
$N_{L^{\theta}}(T^{1})$ acts on $Z$, so $W_{L,\theta}$ acts on
$\mathbb{C}[Z]^{(T)}/\mathbb{C}^{*}$; moreover
$N_{L^{\theta}}(T^{1})$ stabilizes $\overline{T\cdot x_{0}}$, thus
$W_{L,\theta}$ acts on $\mathbb{C}[\overline{T\cdot
x_{0}}]^{(T)}/\mathbb{C}^{*}$. Hence the elements of
$W_{L,\theta}\cdot (C(Z)^{\vee}\cap \chi(S))$ are $T$-weights in
$\mathbb{C}[\overline{T\cdot x_{0}}]^{(T)}/\mathbb{C}^{*}$.   To
prove that $\overline{T\cdot x_{0}}$ is normal it is sufficient to
show that the elements of $W_{L,\theta}\cdot (C(Z)^{\vee}\cap
\chi(S))$ are all the $T$-weights in $\mathbb{C}[\overline{T\cdot
x_{0}}]^{(T)}/\mathbb{C}^{*}$. This holds because of the following
lemma. Notice that $T\cdot x_{0}= S\cdot x_{0}\equiv S$ and
$\mathbb{C}[\overline{T\cdot x_{0}}]^{(T)}/\mathbb{C}^{*}\equiv
\mathbb{C}[\overline{S\cdot x_{0}}]^{(S)}/\mathbb{C}^{*}$. $\square$

\begin{lem}\label{pesiZ}
Let $\lambda$ be a $L$-spherical weight. If there is a
$T$-eigenvector  in $\mathbb{C}[Z]$ of eigenvalue $\lambda$, then
there is a $B_{L}$-eigenvector in $\mathbb{C}[Z]$ of eigenvalue
$\lambda$.
\end{lem}

{\em Proof of Lemma~\ref{pesiZ}.} Let $v$ be a $T$-eigenvector  in
$\mathbb{C}[Z]$ of eigenvalue $\lambda$, we can suppose that it is
contained in an irreducible $L$-subrepresentation $V(\mu)$ of
$\mathbb{C}[Z]$ of highest weight $\mu$. So $\lambda=\mu-\sum
a_{j}\beta_{j}$, where $a_{1},...,a_{n}$ are positive integers and
$\beta_{1},...,\beta_{n}$ are the simple roots of $L$. The cone
$C(Z)$ is generated by simple restricted coroots
$\alpha_{1}^{\vee},...,\alpha_{r}^{\vee}$ and by appropriate vectors
$v_{1},...,v_{m}$ in the antidominant Weyl chamber $-C^{+}$ of
$R_{L,\theta}$. Moreover, for each $i$, $(\lambda,v_{i})=
(\mu,v_{i})+\sum a_{j}(\beta_{j},-v_{i}) \geq (\mu,v_{i})\geq 0$
because $\mu$ belongs to $C(Z)^{\vee}$. On the other hand,
$(\lambda,(\alpha_{i})^{\vee})\geq 0$ for $i=1,...,r$ because
$\lambda$ is dominant, so $\lambda$ belongs to $C(Z)^{\vee}$. Thus
there is an eigenvector in $\mathbb{C}[Z]^{(B_{L})}$ of
$T$-eigenvalue $\lambda$. $\square$

Now, we prove the necessity of the conditions.

\begin{lem}\label{Tx0 smooth} If
$Z$ is smooth, then $[L,L]$ has type $\prod A_{i}$ and
$\overline{T\cdot x_{0}}$  is smooth.
\end{lem}

We begin considering a particular case.

\begin{lem}\label{analisi esplicita rapresentaz}
Suppose that $([L,L],\theta)$ is indecomposable and that $Z(L)$ is
one-dimensional. There are a subgroup $L^{\theta}\subset
H_{L}\subset N_{L}(L^{\theta})$ and an affine smooth embedding of
$L/H_{L}$ if and only if $R_{L}$ has type $A_{l-1}$. If a such
embedding exists, then it is unique, up to an equivariant
isomorphism, and  the closure of $T\cdot x_{0}$ in it is smooth.
\end{lem}

{\em Proof of Lemma~\ref{analisi esplicita rapresentaz}.} Let
$\varsigma$ be the cone in $\chi_{*}(S)\otimes \mathbb{R}$
associated to $\overline{T\cdot x_{0}}$ and let $\varsigma^{\vee}$
be the dual cone. Notice that an affine embedding $Z $ of $L/H_{L}$
is smooth if and only if it is a $L$-representation; in particular
$Z$ must be  an irreducible spherical representation with dimension
lesser than $dim\,[L,L]$. If moreover $[L,L]/[L,L]^{\theta}$ is
isomorphic to a simple group $\dot{L}$, i.e. $[L,L]=\dot{L}\times
\dot{L}$ and $\theta(x,y)=(y,x)$, then $Z$ is isomorphic, as
$[L,L]$-representation,  to $V(\lambda)\otimes
V(-\omega_{0}(\lambda))$. Here $\omega_{0}$ is the longest element
of $W_{\dot{L}}$ and $V(\lambda)$ is a $\dot{L}$-representation with
dimension lesser than $dim\,\dot{L}$. One can easily show that,
given any $L$, there is a such $L$-representation  if and only if
$R_{L,\theta}$ has type $A_{l-1}$ (see   the answer to ex. 24.52 in
\cite{FH}, p. 531 for a list of the irreducible representations of a
simple $L$ with dimension at most $dim\, L$).  In this case $Z$ is
unique, up to twist the $L$-action by a diagram automorphism. More
precisely $Z$ is isomorphic, as $[L,L]$-representation, either to
$V(\omega_{l})$ or to $V(\omega_{1})$. Explicitly we have the
following possibilities for $Z$ (as $[L,L]$-representation): i)
$[L,L] $ is $SL(\mathbb{C}^{n})$, $[L,L]^{\theta}$ is
$SO(\mathbb{C}^{n})$ and $Z$ is either $S^{2}(\mathbb{C}^{n}) $ or
$S^{2}(\mathbb{C}^{n})^{*}\cong S^{2}(\bigwedge^{n-1}\mathbb{C}^{n})
$; ii)   $[L,L] $ is $SL(\mathbb{C}^{2n})$, $[L,L]^{\theta}$ is
$Sp(\mathbb{C}^{2n})$ and  $Z$ is either
$\bigwedge^{2}(\mathbb{C}^{2n})$ or
$\bigwedge^{2}(\mathbb{C}^{2n})^{*}\cong
\bigwedge^{2n-2}(\mathbb{C}^{2n})$; iii) $[L,L] $ is
$SO(\mathbb{C}^{2n+1} )$, $[L,L]^{\theta}$ is $SO(\mathbb{C}^{2n})$
and $Z$ is $\mathbb{C}^{2n+1}$; iv) $[L,L]$ is
$SO(\mathbb{C}^{2n})$, $[L,L]^{\theta}$ is $SO(\mathbb{C}^{2n-1})$
and $Z$ is $\mathbb{C}^{2n}$; v) $[L,L]$ is $E_{6}$,
$[L,L]^{\theta}$ is $F_{4}$ and  $Z$ is either $V(\omega_{1})$ or
$V(\omega_{6})$. Remark that $Z(L)^{0}$ acts by twice a basic
character. Now, we do a case-to-case analysis to prove that, given a
smooth $Z$, the subvariety $\overline{T\cdot x_{0}}$ is smooth. We
denote by $\{f_{1},...,f_{n}\}$ the dual basis of the canonic basis
$\{e_{1},...,e_{n}\}$ of $\mathbb{C}^{n}$.

i) If $([L,L],\theta)$ has type $AI$, let $Q$ be the symmetric
bilinear form $\sum f_{i}^{2}$ over $\mathbb{C}^{n}$. Then the
involution is defined by $\theta(A)=(A^{t})^{-1}$, $L/L^{\theta}$ is
$GL(\mathbb{C}^{n})/O(\mathbb{C}^{n}, Q)$ and $Z$ is, up to
isomorphism, $S^{2}(\mathbb{C}^{n})^{*}$. The subgroup $T\subset
GL(\mathbb{C}^{n})$ of diagonal matrices is a maximally split torus,
the vector fixed by $L^{\theta}$ is $\sum_{i=1}^{n} f_{i}^{2}$ and
the closure of $T\cdot \sum_{i=1}^{n} f_{i}^{2}$ is
$\bigoplus_{i=1}^{n}\mathbb{C} f_{i}^{2}$.

ii) If $(L,L],\theta)$ has type $AII$, let $Q$ be the antisymmetric
bilinear form $\sum_{i=1}^{n}f_{i}\wedge f_{i+n}$ over
$\mathbb{C}^{2n}$. Then $L/L^{\theta}$ is
$GL(\mathbb{C}^{2n})/Sp(\mathbb{C}^{2n}, Q)$ and $Z$ is, up to
isomorphism, $\bigwedge^{2}(\mathbb{C}^{2n})^{*}$. The subgroup $T$
of diagonal matrices is a maximally split torus, the vector fixed by
$L^{\theta}$ is $\sum_{i=1}^{n}f_{i}\wedge f_{i+n}$ and the closure
of $T\cdot (\sum_{i=1}^{n}f_{i}\wedge f_{i+n})$ is the vector space
spanned by $f_{1}\wedge f_{n+1},...,f_{n}\wedge f_{2n}$.

iii) If  $([L,L],\theta)$ has type $BII$, let $Q$ be the symmetric
bilinear form $\sum_{i=1}^{n}f_{i}\wedge f_{i+n}+f_{2n+1}^{2}$ over
$\mathbb{C}^{2n+1}$. Let $\varphi$ be the linear transformation of
$\mathbb{C}^{2n+1}$ such that
\[ \varphi(e_{i})=\begin{cases}
   e_{n+1}&\textrm{if } i=1\\
e_{1}&\textrm{if } i=n+1 \\
e_{i}&\textrm{otherwise}
  \end{cases}.\]
Then $L/L^{\theta}$ is $SO(\mathbb{C}^{2n+1},
Q)\times\mathbb{C}^{*}/ SO((\bigoplus_{i\neq
1,n+1}\mathbb{C}e_{i}\oplus\mathbb{C}(e_{1}+e_{n+1})), Q)\times\{\pm
id\}$ and the involution over $SO(\mathbb{C}^{2n+1}, Q)$ is the
conjugation by the $\varphi$. The representation $Z$ is
$\mathbb{C}^{2n+1}$, over which $\mathbb{C}^{*}$ acts by twice a
basic character. The subgroup $T$ generated by the diagonal matrices
of $SO(\mathbb{C}^{2n+1}, Q)$ and by the connected center is a
maximally split torus, the vector fixed by $L^{\theta}$ is
$e_{1}-e_{n+1}$ and the closure of $T\cdot (e_{1}- e_{n+1})$ is
$\mathbb{C}e_{1}\oplus\mathbb{C}e_{n+1}$.

iv) If $([L,L],\theta)$ has type $DII$, let $Q$ be the  symmetric
bilinear form $\sum_{i=1}^{n }f_{i} f_{i+n}$ over $\mathbb{C}^{2n}$.
Then $L/L^{\theta}$ is $SO(\mathbb{C}^{2n}, Q)\times\mathbb{C}^{*}/
SO((\bigoplus_{i\neq
1,n+1}\mathbb{C}e_{i}\oplus\mathbb{C}(e_{1}+e_{n+1})), Q)\times\{\pm
id\}$ and the subgroup $T$ generated by the diagonal matrices and by
the connected center is a maximally split torus. The representation
$Z$ is $\mathbb{C}^{2n}$, over which $\mathbb{C}^{*}$ acts by  twice
a basic character. The vector fixed by $L^{\theta}$ is
$e_{1}-e_{n+1}$ and the closure of $T\cdot (e_{1}-e_{n+1})$ is the
vector space spanned by $e_{1}$ and $e_{n+1}$.

v) If $([L,L],\theta)$ has type $EIV$, $L/L^{\theta}$ is
$E_{6}\times\mathbb{C}^{*}/F_{4}\times\{\pm id\}$ and  $Z$ is, up to
isomorphism, $V(\omega_{1})$. We will reduce to the already examined
case of $GL_{6}/Sp_{6}$. First, we recall some facts from \cite{LM}.
Let $\mathbb{A}_{\mathbb{C}}$ be either the complexification
$\mathbb{Q}_{\mathbb{C}}$ of the quaternions or the complexification
$\mathbb{O}_{\mathbb{C}}$ of the octonions. Given $a\in \mathbb{A}$,
let $\bar{a}$ be its conjugate. Let
$\mathcal{J}_{3}(\mathbb{A}_{\mathbb{C}})$ be the space of Hermitian
matrices of order three, with coefficients in
$\mathbb{A}_{\mathbb{C}}$:

\[\mathcal{J}_{3}(\mathbb{A}_{\mathbb{C}})=\left\{\left(\begin{matrix}
r_{1}&\bar{x}_{3}&\bar{x}_{2}\\
x_{3}&r_{2}      &\bar{x}_{1}\\
x_{2}&x_{1}      &r_{3}
\end{matrix}\right), r_{i}\in\mathbb{C},x_{i}\in\mathbb{A}_{\mathbb{C}}
\right\}.\]

$\mathcal{J}_{3}(\mathbb{A}_{\mathbb{C}})$  has the structure of a
Jordan algebra with  multiplication  $A\circ B:=\frac{1}{2}(AB+BA)$,
where $AB$ is the usual matrix multiplication. There is a well
defined cubic form on $\mathcal{J}_{3}(\mathbb{A}_{\mathbb{C}})$
called the determinant. The subgroup of
$GL_{\mathbb{C}}(\mathcal{J}_{3}(\mathbb{O}_{\mathbb{C}}))$
preserving the determinant is  $E_{6}$  and contains the subgroup
preserving the Jordan multiplication, namely $F_{4}$. On the other
hand, $SL_{6}$ is contained in $E_{6}$, stabilizes
$\mathcal{J}_{3}(\mathbb{Q}_{\mathbb{C}})$ and its action over
$\mathcal{J}_{3}(\mathbb{Q}_{\mathbb{C}})$ is faithful. The subspace
$\mathcal{J}_{3}(\mathbb{Q}_{\mathbb{C}}) $ of
$\mathcal{J}_{3}(\mathbb{O}_{\mathbb{C}})$ is a (Jordan) subalgebra
and its automorphism group is $Sp_{6}$; moreover one can show that
$Sp_{6}=SL_{6}\cap F_{4}$. There is a simple subgroup $K$ of $F_{4}$
of type $A_{1}$ such that $SL_{6}\cdot K\subset E_{6}$ contains a
maximal torus $T$ of $E_{6}$; $T\times Z(L)^{0}$ is maximally split
because the rank of $E_{6}$ is equal to the sum of the rank of
$F_{4}$ plus the rank of $E_{6}/F_{4}$. The group $F_{4}$ fixes the
identity matrix because $F_{4}$ preserves the Jordan product. Notice
that $ \mathcal{J}_{3}(\mathbb{O}_{\mathbb{C}}) $ is the first
fundamental representation of $E_{6}$, while
$\mathcal{J}_{3}(\mathbb{Q}_{\mathbb{C}}) $ is the second
fundamental representation of $SL_{6}$. Let $T^{1}$ be the maximal
split torus contained in $T\times Z(L)^{0}$, then $T^{1} $ is
contained $SL_{6}\times Z(L)^{0}$, $T^{1}\cdot Id$ is contained
$\mathcal{J}_{3}(\mathbb{Q}_{\mathbb{C}})$ and $Sp_{6}$ fixes $Id$.
(We have assumed that $Z(L)^{0}/\{\pm 1\}$ is the connected center
of $GL_{\mathbb{C}}(\mathcal{J}_{3}(\mathbb{O}_{\mathbb{C}}))$.)
$\square$

{\em Proof of Lemma~\ref{Tx0 smooth}.} We proceed in a similar
manner to the proof of Proposition~\ref{ecc2}. Write $L=\prod
L_{j}\times \widetilde{T}$, where the $L_{j}$ are $\theta$-stable
semisimple normal subgroups of $L$, the $(L_{j},\theta)$ are
indecomposable and $\widetilde{T}$ is a split central torus. Suppose
that there is a smooth affine embedding $Z$ of $L/H_{L}$, for an
opportune subgroup $H_{L}$, and write $Z$ as a sum of irreducible
representations, say $Z=\bigoplus V(\lambda_{h})$. Write $
x_{0}=\sum v_{h}$ with $v_{h}\in V(\lambda_{h})$, so $L\cdot v_{h}$
is dense in $V(\lambda_{h})$.  On the other hand,
$\widetilde{T}\cdot v_{h}$ is contained in $\mathbb{C}^{*}v_{h}$, so
$dim\,L\cdot v_{h}\leq$ $dim\,[L,L]\cdot v_{h}+1$. Write
$\lambda_{i}=\sum\mu_{i,j}+\nu_{i}$, where $\mu_{i,j}$ is a
$L_{j}$-weight  for each $j$ and $\nu_{i}$ is a
$\widetilde{T}$-character;   each  $V(\mu_{i,j})$ contains a $
L_{j}^{\theta}$-fixed vector, in particular $dim V(\mu_{i,j})\geq
dim L_{j}/H_{L_{j}}+1$ for each $i$ and $j$.  The
$[L,L]$-representation $V(\sum_{j=1,...,r}\mu_{i,j})$ is the tensor
product of the $L_{j}$-representations $V( \mu_{i,j})$. Therefore,
if $\mu_{i,j}\neq0$ for at least two index $j_{1}$ and $j_{2}$, then
$dim\, V(\lambda_{i})=dim\, V(\sum\mu_{i,j}) =\prod_{j\,|\,
\mu_{i,j}\neq0} dim\, V(\mu_{i,j})$ \vspace{0.1 mm} $> \sum_{j\,|\,
\mu_{i,j}\neq0}dim\, V(\mu_{i,j})  \geq \sum_{j\,|\, \mu_{i,j}\neq0}
(dim\, (L_{j}/H_{ L_{j}})+1)\geq dim\, L\cdot v_{i}$ (indeed $dim\,
L_{j}/H_{ L_{j}}\geq2$ for each $j$). Thus $V(\lambda_{i})$ is a
smooth affine equivariant embedding of $L_{j_{i}}T_{i}\cdot v_{i}$
for an appropriate index $j_{i}$ and an appropriate one-dimensional
subtorus $T_{i}$ of $\widetilde{T}$; moreover
$Stab_{L_{j_{i}}}(v_{i})\subset N_{L_{j_{i}}}(L_{j_{i}}^{\theta})$.
Every $L_{i}$ acts non-trivially over exactly one $V(\lambda_{i})$
and  $T\cdot x_{0}$ is the product of the varieties $T\cdot v_{i}$,
because $L\cdot x_{0}$ is dense in $Z$. Thus $\overline{T\cdot
x_{0}}$ is the product of the closures of the $T\cdot v_{i}$ (in the
$V(\lambda_{i})$). By the previous lemma such closures are smooth
and the $L_{i}$ have type $A_{l_{i}}$. $\square$

\begin{lem}\label{tx smooth--condizioni}
If    $\overline{Tx_{0}}$  is smooth, then the conditions of the
theorem are verified.
\end{lem}

The hypothesis implies that $\varsigma^{\vee}$ is generated by a
basis $e_{1},...,e_{l}$ of $\chi(S)$. Now, we want to describe
$\mathbb{C}[\overline{T\cdot x_{0}}]$ with relation to the
restricted root system of $(L,\theta)$ (when $\overline{Tx_{0}}$  is
smooth). Recall that the elements of $W_{L,\theta}\cdot
C(Z)^{\vee}\cap \chi(S)$ are the $T$-weights in
$\mathbb{C}[\overline{T\cdot x_{0}}]^{(T)}/\mathbb{C}^{*}$ and that
$\mathbb{C}[\overline{T\cdot x_{0}}]$ is generated by
$S$-seminvariant vectors $v(e_{1}),...,v(e_{l})$  of weight
respectively $e_{1},...,e_{l}$. The symmetric group $S_{l}$ acts on
$\{e_{1},...,e_{l}\}$. Let $\sigma_{i,j}$ be  the transposition that
exchanges $e_{i}$ with $e_{j}$. The action of $W_{L,\theta}$ over
$\chi(S)\otimes\mathbb{R}=\bigoplus \mathbb{R}e_{i}$ is induced by
the action of $N_{L^{\theta}}(T^{1})$ over $\mathbb{C}[T\cdot
x_{0}]^{(T)}$. Moreover $N_{L^{\theta}}(T^{1})$ stabilizes
$\mathbb{C}[\overline{T\cdot x_{0}}]^{(T)}$, so $W_{L,\theta}$
stabilizes the basis $\{e_{1},...,e_{l}\}$ of
$\chi(S)\otimes\mathbb{R}$; in particular $W_{L,\theta}$ exchanges
the $e_{i}$. We want to describe the image of $W_{L,\theta}$ in
$S_{l}$.

\begin{lem}\label{rad+base}
If $\overline{T\cdot x_{0}}$ is smooth, then, for every restricted
$L$-root $\alpha$, there are two indices $i$ and $j$ such that
$\alpha=e_{i}-e_{j}$;  in particular $\sigma_{\alpha}=\sigma_{i,j}$.
Moreover $R_{L,\theta}$ is reduced.
\end{lem}

{\em Proof of the Lemma~\ref{rad+base}.} The orthogonal reflection
$\sigma_{\alpha}$ of $\chi(S)\otimes\mathbb{R}$   stabilizes the
basis $\{e_{1},...,e_{l}\}$. Thus $\sigma_{\alpha}$ can exchange
only two of the $e_{i}$, otherwise the $(-1)$-eigenspace would have
dimension at least two. Suppose now that
$\sigma_{\alpha}=\sigma_{i,j}$, then
$\sigma_{\alpha}(e_{i}-e_{j})=-(e_{i}-e_{j})$, so $e_{i}-e_{j}$ is a
multiple of $\alpha$. Moreover it must be an integral multiple
because $e_{i}-e_{j}=\sigma_{\alpha}(e_{j})-e_{j}$ and $e_{j}$ is a
weight. We can write $\alpha$ as an integral linear combination of
the $e_{j}$ in a unique way, so $e_{i}-e_{j}=\pm\alpha$. Moreover,
there is not a restricted root $\alpha$ such also $2\alpha$ is a
restricted root. Otherwise, by  the previous part of the proof,
there is $i$ and $j$ such that
$\sigma_{\alpha}=\sigma_{2\alpha}=\sigma_{i,j}$ and
$\alpha=2\alpha=\pm(e_{i}-e_{j})$, a contradiction. $\square$

\begin{lem}\label{indicizzazione base}
If $\overline{T\cdot x_{0}}$ is smooth, then, for each fixed $j$,
there is a subset $\{e_{1}^{j},...,e_{l_{j}+1}^{j}\}$ of $\{e_{1}
,...,e_{l }\} $ such that the simple restricted roots in
$R_{L,\theta}^{j}$ are $e_{1}^{j}-e_{2}^{j}$,...,
$e_{l_{j}}^{j}-e_{l_{j}+1}^{j}$; in particular $R_{L,\theta}^{j}$
has type $A_{l_{j}}$.
\end{lem}

{\em Proof of the Lemma~\ref{indicizzazione base}.} The vectors
$e_{i}-e_{j}$ and $e_{h}-e_{k}$ are orthogonal if and only if the
sets $\{i,j\}$ and $\{h,k\}$ are disjoint. Indeed
$\frac{2(e_{h}-e_{k},e_{i}-e_{j})}{(e_{i}-e_{j},e_{i}-e_{j})}(e_{i}-e_{j})
=(e_{h}-e_{k})-\sigma_{e_{i}-e_{j}}(e_{h}-e_{k})$. Let
$\alpha^{j}_{1}$ be an extreme simple root of $R_{L,\theta}^{j}$.
Clearly we can choose $e_{1}^{j}$ and $e_{2}^{j}$ so that
$\alpha^{j}_{1}$ is equal to $e^{j}_{1}-e^{j}_{2}$. Suppose that
there is a subset $\{e_{1}^{j},...,e_{s}^{j}\}$ such that
$\alpha^{j}_{i}$ is equal to $e^{j}_{i}-e^{j}_{i+1}$ for each $i<s$
(with $s<l_{j}+1)$ and let $e_{h}-e_{k}$ be a simple root
$\alpha^{j}_{s}$ which is not orthogonal to the space generated by
$\alpha^{j}_{1}$,...,$\alpha^{j}_{s-1}$. Neither $e_{h}$ nor $e_{k}$
belongs   to $\{e_{1}^{j},...,e_{s-1}^{j}\}$, otherwise the Dynkin
diagram would contain a cycle; thus
$(\alpha^{j}_{s},\alpha^{j}_{i})=0$ for each $i$ strictly lesser
than $s-1$. Hence $(\alpha^{j}_{s},\alpha^{j}_{s-1})\neq0$, so
$e^{j}_{s }$ is equal either to $e_{h}$ or to $e_{k}$. On the other
hand, $e^{j}_{s }$ cannot be equal to $e_{k}$,  otherwise the root
$\sigma_{\alpha^{j}_{s}}(\alpha^{j}_{s-1}) $ is equal to
$e^{j}_{s-1}-e_{h}=e^{j}_{s-1}-e_{k}+e_{k}-e_{h}=
\alpha^{j}_{s-1}-\alpha^{j}_{s}$. Thus $e^{j}_{s }=e_{h}$, we can
define $e^{j}_{s+1 }=e_{k}$ and the first part of the lemma is
proved by induction. In particular,  the Weyl group of
$R_{L,\theta}^{j}$ is isomorphic to $S_{l_{j}+1}$. Therefore
$R_{L,\theta}^{j}$ has type $A_{l_{j}}$, because $R_{L,\theta}$ is
reduced by the Lemma~\ref{rad+base}. $\square$

\emph{Remark 3.} The sets $\{e_{1}^{j},...,e_{l_{j}+1}^{j}\}$ are
pairwise disjoint because the sets $R_{L,\theta}^{j}$ are pairwise
orthogonal. Thus we can index the basis $\{e_{1},...,e_{l}\}$ as
$\{e_{1}^{1},...,e_{l_{1}+1}^{1}, e_{1}^{2},$ $...,$
$e_{l_{2}+1}^{2},$ $...,e_{l_{q}+1}^{q}\}$, so that $\alpha_{i}^{j}=
e_{i}^{j}-e_{i+1}^{j} $ for all $i$ and $j$.

\begin{lem}\label{base C(X)}
Suppose   $\overline{T\cdot x_{0}}$   smooth and let
$\lambda_{i}^{j}$ be $\sum_{h=1}^{i}e_{h}^{j}$ for each $i$ and $j$.
Then $C(Z)^{\vee} \cap\chi(S)$ is generated by the
$\lambda_{i}^{j}$.

\end{lem}

{\em Proof of the Lemma~\ref{base C(X)}.} By the Lemma~\ref{pesiZ}
and the diagram $(2.1)$ we have $C(Z)^{\vee} =
cone(e_{1}^{1},...,e_{l_{q}+1}^{q})\cap C^{+}$. Using the
Lemma~\ref{indicizzazione base},  one can easily show that the
$i$-th fundamental weight $\omega_{i}^{j}$ of $R_{L,\theta}^{j}$ is
equal to
$\frac{1}{l_{j}+1}((l_{j}+1)\lambda_{i}^{j}-i\lambda^{j}_{l_{j+1}})$.
Moreover, $\lambda_{l_{j}+1}^{j}$ is invariant by $W_{L,\theta}$ for
each $j$, so it is orthogonal to span$_{\mathbb{R}}(R_{L,\theta} )$.
Hence $cone(e_{1}^{1},...,e_{l_{q}+1}^{q})\cap C^{+}$ is equal to
$cone(\lambda_{1}^{1},...,\lambda_{l_{q}+1}^{q})$. Indeed the
restriction of $(\omega_{i}^{j},\cdot )$ to
span$_{\mathbb{R}}(R_{L,\theta})$ coincides with the restriction of
$(\lambda_{i}^{j},\cdot )$ for each $i\leq l_{j}$. Thus
$C^{+}\subset cone(\lambda_{1}^{1},...,\lambda_{l_{q}+1}^{q})
+cone(-\lambda_{l_{1}+1}^{1},-\lambda_{l_{2}+1}^{2},...,-\lambda_{l_{q}+1}^{q})$;
moreover, if a vector $\sum a_{i}^{j}\lambda_{i}^{j}=\sum
b_{i}^{j}e_{i}^{j}$ belongs to
$cone(e_{1}^{1},...,e_{l_{q}+1}^{q})$, then
$a_{l_{j}+1}^{j}=b_{l_{j}+1}^{j}$ is positive for each $j$.
$\square$

{\em Conclusion of the proof of the Lemma~\ref{tx
smooth--condizioni}.} The Lemma~\ref{base C(X)} implies the second
condition of   Theorem~\ref{smooth not ecc}.  The type of
$R_{L,\theta}$ is $\prod_{i=1}^{p} A_{l_{i}}$ by the
Lemma~\ref{indicizzazione base}; on the other hand, $rank\, G/H\geq
\sum_{i=1}^{p}(l_{i}+1)$ because of the Lemma~\ref{indicizzazione
base} and of the Remark 3. The last condition of Theorem~\ref{smooth
not ecc} is easily verified using the Lemma~\ref{indicizzazione
base} and the definition of the $\lambda_{i}^{j}$. $\square$

Next, we come to the converse of Lemma~\ref{tx smooth--condizioni}.

\begin{lem}\label{condiz--Tx0 smooth} Suppose that the conditions of the Theorem~\ref{smooth not ecc}
are verified. For all $j$, let $e_{1}^{j}$ be $\lambda_{1}^{j}$ and
let $e_{i+1}^{j}=\sigma_{\alpha_{i}^{j}}e_{i}^{j}$ for every $i\leq
l_{j}$. Then $\varsigma^{\vee}\cap \chi(S)$ is
$\bigoplus_{i,j}\mathbb{Z}^{+}e_{i}^{j}$ and
$e_{i+1}^{j}=\lambda_{i+1}^{j}-\lambda_{i}^{j}$ for each $1\leq
i\leq l_{j}$; in particular $\overline{T\cdot x_{0}}$ is smooth.
\end{lem}

{\em Proof of the Lemma~\ref{condiz--Tx0 smooth}.} Recall that
$\overline{T\cdot x_{0}}$ is the toric variety associated to the
cone $\varsigma^{\vee} =W_{L,\theta}\cdot\,C(Z)^{\vee}$, where
$C(Z)^{\vee}$ is equal to
$cone(\lambda_{1}^{1},...,\lambda_{l_{q}+1}^{q})$; thus
$cone(e_{1}^{1},...,e_{l_{q}+1}^{q})$ is contained in
$\varsigma^{\vee}$. We define $\omega_{0}^{j}:=\lambda_{0}^{j}:=0$
for all $j$ (we use such definition only in the current proof). We
prove that $e_{i+1}^{j}=\lambda_{i+1}^{j}-\lambda_{i}^{j}$ by
induction on $i$. Suppose that this holds for $i-1$ and notice that
$\lambda_{h}^{j}=\omega_{h}^{j}+\frac{h}{l_{j}+1}\lambda_{l_{j}+1}^{j}$
for each $j$ and $h<l_{j}+1$ (here $\omega_{h}^{j}$ is the $h$-th
fundamental weight of $R_{L,\theta}^{j}$). We have
$e_{i+1}^{j}=\sigma_{\alpha_{i}^{j}}(\lambda_{i}^{j}-\lambda_{i-1}^{j})=
\frac{1}{l_{j}+1}\lambda_{l_{j}+1}^{j}+\sigma_{\alpha_{i}^{j}}
(\omega_{i}^{j}-\omega_{i-1}^{j})=
\frac{1}{l_{j}+1}\lambda_{l_{j}+1}^{j}+\omega_{i+1}^{j}-\omega_{i}^{j}=
\lambda_{i+1}^{j}-\lambda_{i}^{j}$, because the
$\lambda_{l_{j}+1}^{j}$ are invariant under $W_{L,\theta}$. The
equality $e_{i+1}^{j}=
\frac{1}{l_{j}+1}\lambda_{l_{j}+1}^{j}+\omega_{i+1}^{j}-\omega_{i}^{j}$
implies that  $\sigma_{\alpha^{j}_{h}}e^{j}_{i+1}\neq e_{i+1}^{j}$
only if $h$ is equal to $i+1$ or to $i$. Thus \vspace{0.1 mm}
$cone(e_{1}^{1},...,e_{l_{q}+1}^{q})$ is stabilized by
$W_{L,\theta}$, hence it is equal to $\varsigma^{\vee}$ because it
contains $C(Z)^{\vee}$. $\square$

\emph{Remark 4.} Suppose that the conditions of the
Theorem~\ref{smooth not ecc} are verified. The previous lemma
together with the Lemma~\ref{indicizzazione base} imply that
$\alpha_{1}^{j}=2\lambda_{1}^{j}-\lambda_{2}^{j}$ and
$\alpha_{i}^{j}=-\lambda_{i-1}^{j}+2\lambda_{i}^{j}-\lambda_{i+1}^{j}$
for each $j$ and $i>1$.

In the   following lemma we  need  that $L$ is as in the Remark 1;
indeed, without such assumption the Lemma~\ref{unicita} is false.

\begin{lem}\label{Tx0 smooth--Z smooth}
If $\overline{T\cdot x_{0}}$ is smooth, then $Z$ is smooth.
\end{lem}

First, we consider a particular case.  Recall that, supposed
$\overline{T\cdot x_{0}}$  smooth, $R_{L,\theta}$ has type $\prod
A_{l_{j}}$.

\begin{lem}\label{unicita}
Suppose that $R_{L,\theta}$ has type $A_{l-1}$ and that $dim\,
Z(L)=1$. Then, there is one subgroup $L^{\theta}\subset H_{L}\subset
N_{L}(L^{\theta})$ and one affine  embedding $Z$ of $L/H_{L}$ such
that the closure of $T\cdot x_{0}$ in $Z$ is smooth (we request that
$H_{L}\cap Z(L)^{0}\equiv\{\pm 1\}\subset\mathbb{C^{*}}\equiv
Z(L)^{0}$). Moreover: i) $H_{L}$ and $Z$ are unique up to an
equivariant isomorphism; ii) $Z$ is smooth.
\end{lem}

{\em Proof of the Lemma~\ref{unicita}.} We have already show in the
proof of the Lemma~\ref{analisi esplicita rapresentaz} that, there
are one subgroup $H'$ and one affine $L$-variety $Z'$ such that: i)
$L^{\theta}\subset H' \subset N_{L}(L^{\theta})$; ii) $Z'$ is an
affine smooth embedding of $L/H'$; such variety is unique up to an
equivariant isomorphism and  the closure of $T\cdot x_{0}$ in it is
smooth. Therefore, by the Lemma~\ref{condiz--Tx0 smooth}, it is
sufficient to prove that there is one $L^{\theta}\subset H_{L}
\subset N_{L}(L^{\theta})$ and one affine embedding of $L/H_{L}$ (up
to equivariant isomorphism) satisfying the conditions of the
Theorem~\ref{smooth not ecc}.

Let $H_{L}$ be a subgroup of $L$ and let $Z$ be an affine embedding
of $L/H_{L}$  satisfying the conditions of the Theorem~\ref{smooth
not ecc}; write $C(Z)=cone((\alpha^{1}_{1})^{\vee},...,$
$(\alpha^{1}_{l-1})^{\vee},v_{1})$. To determine $H_{L}$ is
sufficient to determine $H_{T}$ (see \cite{V1}, \S 2.2, Lemme 2).
First, we determine the subgroup $H_{T\cap [L,L]}$. Observe that
$\bigoplus_{i=1}^{l-1}\mathbb{Z}(\alpha^{1}_{i})^{\vee}=\chi_{*}(T\cap[L,L]/T\cap[L,L]^{\theta})$
has finite index in $\chi_{*}(T\cap[L,L]/H_{T\cap[L,L]})$ and that
$\chi_{*}(T\cap[L,L]/H_{T\cap[L,L]})$ is contained in $\chi_{*}(S)$.
Hence $\chi_{*}(T\cap[L,L]/H_{T\cap[L,L]})$ is
$\bigoplus_{i=1}^{l-1}\mathbb{Z}(\alpha^{1}_{i})^{\vee}$ and
$H_{T\cap[L,L]}$ is $T\cap [L,L]^{\theta}$.  On the other hand we
have supposed that $H_{Z(L)^{0}}=(Z(L)^{0})^{\theta}$, thus
$\chi:=\chi (T\cap[L,L]/H_{T\cap[L,L]}\times Z(L)^{0}/H_{Z(L)^{0}})$
is $\chi (T/T^{\theta})=\bigoplus_{i=1}^{l-1}
\mathbb{Z}\omega^{1}_{i}\oplus \mathbb{Z}\frac{1}{h}\lambda^{1}_{l}$
for an appropriate integer $h$, so $\chi(S)$ has index at least $l$
in $\chi$ (observe that $\chi
(Z(L)^{0})\otimes\mathbb{R}=\mathbb{R}\lambda^{1}_{l}$ because
$\lambda^{1}_{l}$ is orthogonal to all the simple restricted
coroots). There is exactly one subgroup $K$ of
$N_{L}(L^{\theta})/L^{\theta}$ of order at least $l$  not
intersecting neither $N_{[L,L]}([L,L]^{\theta})/[L,L]^{\theta}$ nor
$Z(L)^{0}/\{\pm id\}$. Indeed
$N_{[L,L]}([L,L]^{\theta})/[L,L]^{\theta}$ is isomorphic to the
fundamental group of $R_{L,\theta}$, namely the cyclic group of
order $l$; on the other hand all  the  finite subgroups of
$\mathbb{C}^{*}$ are cyclic. Therefore $H_{L}/L^{\theta}$ is equal
to $K$, so $\chi_{*}(S)$ and $H_{L}$ are univocally determined. By
the Remark 4, we can suppose that
$(\alpha^{1}_{i},v_{1})=-\delta_{i,l-1}$ up to isomorphism. On the
other hand, $\{(\alpha^{1}_{1})^{\vee},...,$
$(\alpha^{1}_{l-1})^{\vee},v_{1}\}$ is a basis of $\chi_{*}(S)$,
thus the image of $v_{1}$ by the projection
$(\chi_{*}(Z(L)^{0})\otimes\mathbb{R})\oplus\bigoplus\mathbb{R}(\alpha_{i}^{1})^{\vee}
\rightarrow \chi_{*}(Z(L)^{0})\otimes\mathbb{R}$ is determined up
the sign. Therefore $C(Z)$ is univocally determined up the
automorphism of $L$ which is the identity on $[L,L]$ and coincides
with $\theta$ on $Z(L)^{0}$. On the other hand, $D(Z)=D(L/H_{L})$
because $Z$ is affine and $dim\,C(Z)=ran\!k\, \chi_{*}(S)$ (see
\cite{Br3}, Corollaire 2, page 51). Therefore,   $Z$ is univocally
determined up to automorphism. $\square$

{\em Proof of the Lemma~\ref{Tx0 smooth--Z smooth}.} Let
$\{(\alpha_{1}^{1})^{\vee},...,(\alpha_{l_{1}}^{1})^{\vee},v_{1},
(\alpha_{1}^{2})^{\vee},$ $...,$
$(\alpha_{l_{2}}^{2})^{\vee},v_{2},$ $...,v_{q}\}$ be the dual basis
of $\{\lambda_{1}^{1},...,\lambda_{l_{1}+1}^{1}, \lambda_{1}^{2},$
$...,$ $\lambda_{l_{2}+1}^{2},$ $...,\lambda_{l_{q}+1}^{q}\}$. Write
$L=\prod L_{i}\times \widetilde{T}$ with $(L_{i},\theta)$
indecomposable and $\widetilde{T}=Z(L)^{0}$. Let $T_{i}$ be the
inverse image in $T^{1}$ of $Im(v_{i})\subset S$ and  define
$L'_{i}=L_{i}T_{i}$. Observe that $\chi_{*}(Z(L)^{0})$ is
$\frac{1}{2}\chi_{*}(Z(L)^{0}/H_{Z(L)^{0}})$, $\chi_{*}(T^{1}\cap
L'_{i}\,/\,H_{T^{1}\cap L'_{i}})$ is $\mathbb{Z}v_{i}\oplus
\bigoplus_{j=1}^{l_{j}+1} \mathbb{Z}(\alpha^{i}_{j})^{\vee}$ and
$\chi_{*}(S)$ is $\bigoplus_{i}\, \chi_{*}(T^{1}\cap
L'_{i}\,/\,H_{T^{1}\cap L'_{i}})$. Thus $L$ is the direct product of
the $L'_{i}$ and $H_{L}=\prod H_{L'_{i}}$. Hence $L/H_{L} =\prod
L'_{i}/H_{L'_{i}}$. We have $C(Z)^{\vee}=\sum_{i=1}^{q}
cone(\lambda^{i}_{1},...,\lambda^{i}_{l_{i}+1})$, hence  $Z$ is the
product of the affine embeddings $Z_{i}$ of $L'_{i}/H_{ L'_{i}}$
corresponding respectively to the colored cones
$(cone(\lambda^{i}_{1},...,\lambda^{i}_{l_{i}+1}))^{\vee},\rho^{-1}
(\{\alpha^{\vee}_{1},...,\alpha^{\vee}_{l_{i}}\}))$. These
embeddings are smooth by the previous lemma, so $Z$ is smooth.
$\square$

{\em Remark 5.} Let $X$ be a smooth simple symmetric  variety and
suppose for simplicity that $\chi(S)$ has rank equal to the rank of
$R_{L}$ plus one, in particular that $R_{L}$ is irreducible. We want
to remark that there is only one way to index the basis of $\chi(S)$
that generates $C(X)^{\vee}$ so that it verifies the third
condition. First, observe that  indexing  this basis is equivalent
to indexing the dual basis of $\chi_{*}(S)$. Write $C(X)\cap
\chi_{*}(S)=\bigoplus\mathbb{Z}^{+}(\alpha_{i}^{1})^{\vee}\oplus
\mathbb{Z}^{+}v_{1}$. If we request that $(\alpha^{1}
_{i},\alpha^{1}_{j} )\neq 0$  only if $i\in\{j-1,j,j+1\}$, then the
index of $\{(\alpha_{1}^{1})^{\vee},
...,(\alpha_{l-1}^{1})^{\vee}\}$ is defined up the non trivial
automorphism of the Dynkin diagram. The third condition of the
theorem implies that
$\alpha^{1}_{1}=2\lambda^{1}_{1}-\lambda^{1}_{2}$ and
$\alpha^{1}_{i}=-\lambda^{1}_{i-1}+2\lambda^{1}_{i}-\lambda^{1}_{i+1}$
for each $i>1$; in particular
$(\alpha^{1}_{i},v_{1})=-\delta_{l-1,i}$. Hence the condition
$(\alpha^{1}_{l-1},v_{1})=-1$ let us  determine completely the
indexing of the basis $\{(\alpha_{1}^{1})^{\vee},
...,(\alpha_{l-1}^{1})^{\vee},v\}$.

\section{Smooth Complete Symmetric
Varieties with Picard Number  One}

\emph{In this section we always suppose $G$ semisimple.} We classify
the smooth complete symmetric varieties with Picard number one and
we prove that they are all projective. Let $X$ be a smooth complete
symmetric variety of rank $l$.   Let $r$ be the number of colors
which do not belong to $X$ and let $m$ be the number of
one-dimensional cones which are faces of a colored cone belonging to
the colored fan of $X$. Observe that we consider also the
one-dimensional faces which do not intersect the valuation cone. It
is known that the Picard group of $X$ is free and its rank is equal
to $r+m-l$. In fact, the rank of $Pic(X)$ is equal to cardinality of
$D(G/H)\backslash\, D(X)$ minus the rank of $G/H$ plus the rank of
the group composed by the functions over $C(X)$ which are linear
over each colored cone of $X$ and assume integral value over
$\chi_{*}(S)\cap C(X)$ (see \cite{Br3}, \S 5.2, Th\'{e}or\`{e}me);
moreover the maximal colored cones are generated by a basis of
$\chi_{*}(S)$. If $Pic(X)$ has rank one, then there is at most one
color which does not belong to $X$. We have two cases: i) there is
one color which does not belong to $X$ and $X$ is simple ($r=1$ and
$m=l$); ii) all the colors belong to $X$ and $m=l+1$.

\begin{lem}\label{l-colors}
If $G/H$  has a smooth completion $X$ with Picard number one, then
$G/H$ has exactly $l$ colors. Moreover, any maximal simple open
$G$-stable subvariety of $X$ has still Picard number one.
\end{lem}

{\em Proof.} There are exactly $l+1$ $B$-stable prime divisors on
$X$ because $r+m=l+1$  and $\rho$ is injective over $D(X)$. Among
them, there must be at least one $G$-stable divisor.  (The cones
spanned only by colors intersect the valuation cone just  in 0). On
the other hand, there are at least $l$ colors. Hence there are
exactly $l$ colors corresponding to $l$ distinct restricted coroots.
$\square$

Notice that $N(X)$ contains one vector, say $v$. Let $\widetilde{X}$
be a maximal simple open $G$-subvariety of $X$.  In the following we
denote by  $L$ the standard Levi subgroup of the stabilizer of
$\widetilde{X}_{B}$. First of all, we prove that either
$R_{G,\theta}$ is irreducible or has type $A_{1}\times A_{1}$; next
we will do a case-by-case analysis. In some cases we will classify
first the possibilities for $\widetilde{X}$ and in a second time the
possibilities for $X$. Observe that in the previous section $X$ was
always   a simple variety.

\begin{lem}\label{L-irr}\begin{enumerate}
\item[(i)] The root system $R_{L,\theta}$ is irreducible
and has type $A_{l-1}$.
\item[(ii)] Denote by $\alpha$  the simple restricted root such that
$\alpha^{\vee}$ does not belong to $\rho(D(\widetilde{X}))$. Then
$\alpha$ must be an endpoint of the Dynkin diagram of
$R_{G,\theta}$.
\end{enumerate}
\end{lem}

{\em Proof.} The rank of $R_{G,\theta}$ is equal to rank of
$R_{L,\theta}$ plus one, \vspace{0.2 mm} because
$\rho(D(\widetilde{X}))$ contains $l-1$ elements and $\rho$ is
injective. On the other hand, as $\widetilde{X}$ is smooth, the
$rank\,(R_{G,\theta})$ is also greater (or equal) than
$rank\,(R_{L,\theta})$ plus the number of connected components of
the Dynkin diagram of $R_{L,\theta}$. Thus the first point of the
lemma follows. To prove the second point observe that the Dynkin
diagram of $R_{L,\theta}$ is obtained from the one  of
$R_{G,\theta}$ by removing the vertex corresponding to $\alpha$.
$\square$

\begin{lem}\label{G-irr}
The root system $R_{G,\theta}$ is irreducible or its type is
$A_{1}\times A_{1}$.
\end{lem}

{\em Proof.} Suppose  that  $R_{G,\theta}$ is reducible, then its
type is either $A_{l-1}\times A_{1}$ or $A_{l-1}\times BC_{1}$,
because of Lemma~\ref{L-irr}. We index the simple roots of
$R_{G,\theta}$ so that $\alpha_{l}$ is orthogonal to the other
simple roots and we denote by $\omega_{1}$ the first fundamental
weight of the subroot system of type $A_{l-1}$; in particular we
assume $D(\widetilde{X})=\{D_{\alpha_{i}^{\vee}}\}_{i<l}$. The
variety $\widetilde{X}$ cannot be complete, otherwise
$C(\widetilde{X})$ contains $-\omega^{\vee}_{1}$. Notice that
$C(\widetilde{X})$ contains also $\omega^{\vee}_{1}$, because
$\omega^{\vee}_{1}$ is a linear combination of
$\alpha^{\vee}_{1},...,\alpha^{\vee}_{l-1}$ with positive
coefficients. There is another maximal simple subvariety with colors
$\{D_{\alpha_{j}^{\vee}}\}_{j\neq i}$ for an appropriate $i<l$,
because $\widetilde{X}$ is not complete and $X$ contains exactly one
$G$-stable divisor. Hence, by the Lemma~\ref{L-irr}, there is $i<l$
such that $\alpha_{1},...,\widehat{\alpha}_{i},...,\alpha_{l}$
generate a root system of type $A_{l-1}$. This is possible only if
$l$ is equal to two and $R_{G,\theta}$ has type  $A_{1}\times
A_{1}$. $\square$

{\em Remarks 6:} i) by the previous two lemmas, the variety $X$ may
be non-simple only if $R_{G,\theta}$ has type $A_{1}\times A_{1}$,
$A_{l}$, $B_{2}$, $D_{l}$ or $G_{2}$; ii) because of the Remark 5,
if
$\rho(D(\widetilde{X}))=\{\alpha_{1}^{\vee},...,\alpha_{l-1}^{\vee}\}$,
then $(\alpha_{i}^{\vee},v)$ is equal either to $-\delta_{i,1}$ or
to $-\delta_{i,l-1}$ for each $i\leq l-1$; iii) if $G/H$ is
Hermitian then, by the Lemma~\ref{l-colors}, $H=N_{G} (G^{\theta}) $
and there are no  exceptional roots.

In the following we do a case-by-case analysis.

i) Suppose  that $R_{G,\theta}$ is reducible. We have seen that it
has type $A_{1}\times A_{1}$ and that $X$ is not simple. Thus $X$ is
covered by two simple open subvarieties whose associated colored
cones are respectively
$(cone(\alpha^{\vee}_{1},v),\{D_{\alpha^{\vee}_{1}}\})$ and
$(cone(\alpha^{\vee}_{2},v),\{D_{\alpha^{\vee}_{2}}\})$; so
$v=-\omega_{1}^{\vee}-\omega_{2}^{\vee}$ because of the  Remarks 6.
Observe that the lattice generated by $\alpha^{\vee}_{1}$ and $v$ is
equal to the lattice generated by $\alpha^{\vee}_{2}$ and $v$; more
precisely, it is the lattice
$\mathbb{Z}2\omega_{1}^{\vee}\oplus\mathbb{Z}(\omega_{1}^{\vee}+\omega_{2}^{\vee})$.
Let $\{\lambda_{1},\lambda_{2}\}$ be the dual basis of
$\{\alpha^{\vee}_{1},v\}$. We have
$\lambda_{1}=\frac{1}{2}\alpha_{1}-\frac{1}{2}\alpha_{2}$,
$\lambda_{2}=-\alpha_{2}$ and
$\frac{1}{2}(2\lambda_{1}-\lambda_{2})=\frac{1}{2}\alpha_{1}=\omega_{1}$.
We can proceed in the same manner for the dual   basis of
$\{\alpha^{\vee}_{2},v\}$. We have proved that in this case $X$ is
smooth.

\emph{In the following we always assume  $R_{G,\theta}$
irreducible.} Moreover, we number the simple restricted roots as in
\cite{H}.
\begin{lem} \label{unic reticolo} Let  $(G,\theta)$  be
an indecomposable  involution and  let
$(cone(\alpha^{\vee}_{1},...,$ $\alpha^{\vee}_{l-1},$ $v),
\rho^{-1}(\{\alpha^{\vee}_{1},...,\alpha^{\vee}_{l-1}\}))$ be a
colored cone in $\chi_{*}(T^{1})\otimes\mathbb{R}$. There is at most
one subgroup $G^{\theta}\subset H\subset N_{G}(G^{\theta})$ such
that the embedding $\widetilde{X}$ of $G/H$ corresponding to the
previous colored cone is smooth.
\end{lem}

{\em Proof of Lemma~\ref{unic reticolo}.} Suppose that a such group
exists; in particular $R_{L,\theta}$ has type $A_{l-1}$. We can
suppose that $\{\alpha^{\vee}_{1},...,\alpha^{\vee}_{l-1},v\}$ is a
basis of $\chi_{*}(S)$; let $\{\lambda_{1},...,\lambda_{l}\}$ be the
dual basis. We have to prove that $\chi_{*}(S)$, or equivalently
$v$, are univocally determined. Notice that $v$ is determined up to
a multiplicative constant because it generates a one-dimensional
face of the fixed colored cone. Moreover, if we substitute $v$ with
a multiple, $\lambda_{1},...,\lambda_{l-1}$ do not change because
$(\lambda_{i},v)=0$ for each $i<l$. Thus, it is sufficient to prove
that there is a unique possibility for $\lambda_{l}$. By the
Theorem~\ref{smooth not ecc}, $\lambda_{l}$ is equal to
$l(\lambda_{1}-\widetilde{\omega}_{1})$, thus it is univocally
determined (here $\widetilde{\omega}_{1}$ is the first fundamental
weight of the root system generated by
$\alpha_{1},...,\alpha_{l-1}$). $\square$

ii) Suppose  that $R_{G,\theta}$ has rank one. Then, for each $H$
such that $G^{\theta}\subset H\subset N_{G}(G^{\theta})$, there is a
unique (non trivial) embedding and it is simple, projective, smooth,
with Picard number at most two and $D(X)$ equal to the empty set.
The group $Pic(X)$ has rank two if and only if $G/H$ is Hermitian
and $H=G^{\theta}$. In this case we have two possibilities: i)
$(G,\theta)$ has type $AI$,  $G/H$ is isomorphic to $SL_{2}/SO_{2}$,
the restricted root system has type $A_{1}$ and $X$ is isomorphic to
$\mathbb{P}^{1}\times \mathbb{P}^{1}$; ii) $(G,\theta)$ has type
$AIV$, $G/H$ is isomorphic to $SL_{n+1}/S(L_{1}\times L_{n})$, the
restricted root system has type $BC_{1}$, $X$ is exceptional and
isomorphic to $\mathbb{P}^{n}\times (\mathbb{P}^{n})^{*}$.

\emph{In the following we always assume that $R_{G,\theta}$  has
rank at least two.} We will work with the standard inner product
over $\mathbb{R}^{n}$  and we will denote by $\{e_{1},...,e_{n}\}$
the usual orthonormal basis. One can realize a root system of type
$A_{l-1}$ as the set $\{ e_{i}- e_{j}\,\mid\, i\neq j\}$ contained
in the vector space $\{\sum a_{i}e_{i}\in \mathbb{R}^{l}\,|\,\sum
a_{i}=0\}$. Moreover, $\{e_{1}-e_{2},...,e_{l-1}-e_{l}\}$ is a basis
of the root system and the  $j$-th fundamental weight is
$\frac{1}{l}(l\sum_{i=1}^{j}e_{i}-j\sum_{i=1}^{l}e_{i})$.

iii) Suppose  that $R_{G,\theta}$ has type $A_{l}$. As before, we
can realize $R_{G,\theta}$ as the set $\{ e_{i}- e_{j}\,\mid\, i\neq
j\}\subset \mathbb{R}^{l+1}$; in particular we can identify the
restricted roots with the corresponding coroots. Up to an
automorphism of the Dynkin diagram, we can suppose that
$\rho(D(\widetilde{X}))=\{\alpha_{1},...,\alpha_{l-1}\}$. First, we
consider the case where $X$ is simple, i.e. $\widetilde{X}$ is
complete and coincides with $X$. In this case $C(\widetilde{X})$ is
generated by $-\omega_{1} $ and $\rho(D(\widetilde{X}))$. Indeed
$-\omega_{2}= \alpha_{1}-2\omega_{1}$ and
$-\omega_{j+1}=\sum_{i=1}^{j}\alpha_{j}-\omega_{1}-\omega_{j}$ for
each $j=2,...,l-1$. We set $f_{i}=-e_{l-i+1}$ and $
\gamma_{i}=\alpha_{l-i}$, so that $\gamma_{i}=f_{i}-f_{i+1}$ and
$-\omega_{1} =\frac{1}{l+1}( (l+1)f_{l}-\sum_{i=0}^{l}f_{i})$ (we
allow $i$ to be zero). We have renamed the simple restricted roots
so that $(\gamma_{l-1},-\omega_{1})=-1$. Let
$\{\lambda_{1},...,\lambda_{l}\}$ be the dual basis of
$\{\gamma_{1},...,\gamma_{l-1},-\omega_{1}\}$. We have $\lambda_{j}=
(\sum_{i=1}^{j} f_{i})-jf_{0} $, so
$(1/l)(l\lambda_{j}-j\lambda_{l})
=(1/l)(l\sum_{i=1}^{j}f_{i}-j\sum_{i=1}^{l}f_{i})$ is the $j$-th
fundamental weight $\widetilde{\omega}_{j}$ of the \vspace{0.45 mm}
root system generated by $\gamma_{1},...,\gamma_{l-1}$. Hence, the
symmetric variety with
$\rho(D(\widetilde{X}))=\{\alpha_{1},...,\alpha_{l-1}\}$,
$N(\widetilde{X})=\{-\omega_{1}\}$ and $\chi_{\ast}(S)$ equal to the
weight lattice of $R^{\vee}_{G,\theta}$ is smooth.  It is the only
possibility for a simple $X$  because of the Lemma~\ref{unic
reticolo}.

Now we  classify the $X$ which are not simple. In this case $X$ has
two closed orbits; indeed the maximal colored cones of the colored
fan of $X$ are
$(cone(\alpha_{1},...,\alpha_{l-1},v),\{D_{\alpha^{\vee}_{1}},...,$
$D_{\alpha^{\vee}_{l-1}}\})$ and
$(cone(\alpha_{2},...,\alpha_{l},v),\{D_{\alpha^{\vee}_{2}},...,$
$D_{\alpha^{\vee}_{l}}\})$. Let $b\omega_{1}$ be the primitive
vector of $\mathbb{R}^{+}\omega_{1}$; we have \vspace{0.1 mm} $1\leq
b\leq l+1$. Suppose first that $(v,\alpha_{l-1})=-\delta_{i, 1}$ for
each $i<l$; then $v=\sum_{i=1}^{l-1} a_{i}\gamma_{i}-b\omega_{1}$,
where the $\gamma_{i}$ are as before. The $a_{i}$ are positive
integers, because $-C^{+}\cap
\chi_{*}(S)\subset\mathbb{Z}^{+}(-\omega_{1})+\sum\mathbb{Z}^{+}\gamma_{i}$;
the coefficient of $v$ with respect to $-\omega_{1}$ is $b $ because
$-b\omega_{1}$ belongs \vspace{0.1 mm} to the the lattice generated
by $\gamma_{1},...,\gamma_{l-1}$ and $v$. Let
$\{\widetilde{\lambda}_{1},...,\widetilde{\lambda}_{l}\}$ be the
dual basis \vspace{0.1 mm} of $\{\gamma_{1},...,\gamma_{l-1},v\}$;
we have $\widetilde{\lambda}_{l}=(1/b)\lambda_{l}$ and
$\widetilde{\lambda}_{i}=\lambda_{i}-(a_{i}/b)\lambda_{l}$ for each
$i<l$. Because $\widetilde{X}$ is smooth we have \vspace{0.3 mm}
$\widetilde{\omega}_{j}=\frac{1}{l}
(l\widetilde{\lambda}_{j}-j\widetilde{\lambda}_{l})=
\frac{1}{l}(l\lambda_{j}-j\lambda_{l})+\frac{1}{l}(j-\frac{la_{j}}{b}-\frac{j}{b})\lambda_{l}$,
so $a_{j}=\frac{(b-1)j}{l}$. In particular $l$ divides $b-1$, so
either $b=1$ or $b=l+1$.   If $b=1$, then the $a_{j}$ are all zero
and $\widetilde{X}$ is complete, a contradiction. If $b=l+1$, then
$a_{j}=j$ for each $j$ and $v=-\omega_{1}-\omega_{l}$ is a root.
Thus $\chi(S)$ must be the weight lattice of $R_{G,\theta}$. Notice
that $-\omega_{1}-\omega_{l}=-\alpha_{1}-...-\alpha_{l}$, so
$\{\alpha_{l-1},...,\alpha_{1},-\omega_{1}-\omega_{l}\}$ is a basis
of $\chi_{*}(S)$. Observe that $-\omega_{1}-\omega_{l}$ is fixed by
the non trivial automorphism of the Dynkin diagram. Next, suppose
that $(\alpha_{i},v)=-\delta_{i,l-1}$ for each $i<l$, so
$v=-\omega_{l-1}-a\omega_{l}$. We exclude this case because the
simple variety associated to $cone(\alpha_{2},...,\alpha_{l},v)$ is
not   smooth.

We have proved that, if $X$ is  not simple, then it is covered by
two simple varieties   corresponding to the colored cones
$(cone(\alpha_{1},...,\alpha_{l-1},$
$-\omega_{1}-\omega_{l}),\{D_{\alpha^{\vee}_{1}},...,$
$D_{\alpha^{\vee}_{l-1}}\})$ and
$(cone(\alpha_{2},...,\alpha_{l},-\omega_{1}-\omega_{l}),$
$\{D_{\alpha^{\vee}_{2}},...,$ $D_{\alpha^{\vee}_{l}}\})$. If $l$ is
at least three, then $-\omega_{2}$ does not belongs neither to
$cone(\alpha_{1},...,\alpha_{l-1},-\omega_{1}-\omega_{l})$ nor  to
$cone(\alpha_{2},...,\alpha_{l},-\omega_{1}-\omega_{l})$; thus $l$
must be two. In this case,  the previous colored cones define a
complete variety because $cone(-\alpha_{1},-\alpha_{1}-
\alpha_{2})\subset cone(\alpha_{2},-\omega_{1}- \omega_{2})$.

iv) Suppose   that $R_{G,\theta}$ has type $B_{l}$ with $l>2$; then
$X$ must be simple. Observe that the dual root system has type
$C_{l}$ and
$\rho(D(X))=\{\alpha^{\vee}_{1},...,\alpha^{\vee}_{l-1}\}$. We can
realize $R_{G,\theta}$ as the set $\{\pm(e_{i}\pm e_{j})\,|\, i\neq
j\}\cup\{\pm e_{i}\}\subset \mathbb{R}^{l}$; moreover we can suppose
that the basis of $R_{G,\theta}$ is
$\{e_{1}-e_{2}$,...,$e_{l-1}-e_{l},e_{l}\}$. The dual root system is
$\{\pm(e_{i}\pm e_{j})\,|\, i\neq j\}\cup\{\pm 2e_{i}\}\subset
\mathbb{R}^{l}$ and has basis
$\{e_{1}-e_{2}$,...,$e_{l-1}-e_{l},2e_{l}\}$. Notice that the cone
generated by $\rho(D(X))$ and $-C^{+}$ is equal to the cone
$cone(\alpha^{\vee}_{1},...,\alpha^{\vee}_{l-1},-\omega^{\vee}_{1})$.
Indeed $-\omega^{\vee}_{2}= \alpha^{\vee}_{1}-2\omega^{\vee}_{1}$
and $-\omega^{\vee}_{i}=
\sum_{j=1}^{i-1}\alpha^{\vee}_{j}-\omega^{\vee}_{1}-\omega^{\vee}_{i-1}$
for $i\geq3$. If  $\chi(S)$ is the root lattice of $R_{G,\theta}$,
then $C(X)\cap\chi_{*}(S)$ is generated by
$\{\alpha^{\vee}_{1},...,\alpha^{\vee}_{l-1},-\omega^{\vee}_{1}\}$,
which is a basis of $\chi_{*}(S)$. Now we can proceed as in the case
of $A_{l}$, because the Cartan numbers
$(\alpha_{i},\alpha^{\vee}_{j})$ with $j<l$ are equal to  the ones
of $A_{l}$; explicitly  we define  $f_{i}=-e_{l-i+1}$ and
$\gamma_{i}=\alpha_{l-i}$, so that
$-\omega_{1}=-\omega_{1}^{\vee}=f_{l}$, $\gamma_{0}=-f_{1}$ and
$\gamma_{i}=\gamma_{i}^{\vee}=f_{i}-f_{i+1}$ for each $i=1,...,l-1$.
Let $\{\lambda_{1},...,\lambda_{l}\}$ be the dual basis of
$\{\gamma^{\vee}_{1},...,\gamma^{\vee}_{l-1},-\omega^{\vee}_{1}\}$;
we have $\lambda_{j}= \sum_{h=1}^{j }f_{h}$, so
$(1/l)(l\lambda_{j}-j\lambda_{l}) =
(1/l)(l\sum_{h=1}^{j}f_{h}-j\sum_{h=1}^{l }f_{h})$ is the $j$-th
fundamental weight  of the root system generated by
$\gamma_{1},...,\gamma_{l-1}$. We have proved that in this case $X$
is smooth. It is the only possibility because of Lemma~\ref{unic
reticolo}.

v) Suppose now that $R_{G,\theta}$ has type $C_{l}$. Then $X$ must
be simple and $\rho(D(X))$ must be equal to
$\{\alpha^{\vee}_{1},...,\alpha^{\vee}_{l-1}\}$. Observe that the
dual root system has type $B_{l}$. We can realize $R_{G,\theta}$ as
$\{\pm(e_{i}\pm e_{j})\,|\, i\neq j\}\cup\{\pm2e_{i}\}\subset
\mathbb{R}^{l}$ with basis
$\{e_{1}-e_{2}$,...,$e_{l-1}-e_{l},2e_{l}\}$; the dual root system
is $\{\pm(e_{i}\pm e_{j})\,|\, i\neq j\}\cup\{\pm e_{i}\}$. The cone
generated by $\rho(D(X))$ and $-C^{+}$ is
$cone(\alpha^{\vee}_{1},...,\alpha^{\vee}_{l-1},-\omega^{\vee}_{1})$.
Indeed $-\omega^{\vee}_{2}= \alpha^{\vee}_{1}-2\omega^{\vee}_{1}$,
$-\omega^{\vee}_{i}=
\sum_{j=1}^{i-1}\alpha^{\vee}_{j}-\omega^{\vee}_{1}-\omega^{\vee}_{i-1}$
for $2<i<l$ and $-2\omega^{\vee}_{l}=
\sum_{j=1}^{l-1}\alpha^{\vee}_{j}-\omega^{\vee}_{1}-\omega^{\vee}_{l-1}$.
The lattice $\chi_{*}(S)$ is generated by the restricted coroots
because $\omega^{\vee}_{l}$ does not belong to
$\mathbb{Z}(-\omega^{\vee}_{1})\oplus\bigoplus_{i=1}^{l-1}\mathbb{Z}\alpha^{\vee}_{i}$;
thus $G=G^{\theta}$. We have to exclude the cases where $G/H$ is
Hermitian, because in such cases  $X$ would have rank two. The set
$\{\alpha^{\vee}_{1},...,\alpha^{\vee}_{l-1},-\omega^{\vee}_{1}\}$
is a basis of $\chi_{*}(S)$ because
$-\omega^{\vee}_{1}=-\sum_{i=1}^{l}\alpha^{\vee}_{i}$. We define
$f_{i}=-e_{l-i+1}$ and $\gamma_{i}=\alpha_{l-i}$, so that
$-\omega_{1}=-\omega_{1}^{\vee}=f_{l}$, $\gamma_{0}=-2f_{1}$ and
$\gamma_{i}=\gamma_{i}^{\vee}=f_{i}-f_{i+1}$ for each $i=1,...,l-1$.
Let $\{\lambda_{1},...,\lambda_{l}\}$ be the dual basis of
$\{\gamma^{\vee}_{1},...,\gamma^{\vee}_{l-1},-\omega_{1}\}$, then
$\lambda_{j}= \sum_{h=1}^{j }f_{h}$ and
$(1/l)(l\lambda_{j}-j\lambda_{l}) =(1/l)(l\sum_{i=1}^{j
}f_{i}-j\sum_{i=1}^{l }f_{i})$ is the $j$-th fundamental weight of
the root system generated by $\gamma_{1},...,\gamma_{l-1}$.  Thus,
$X$ is smooth if
$C(X)=cone(\gamma^{\vee}_{1},...,\gamma^{\vee}_{l-1},-\omega_{1}\}$
and $\chi_{*}(S)=\bigoplus_{i=0}^{l}\mathbb{Z}\gamma^{\vee}_{i}$; it
is the only possibility because of Lemma~\ref{unic reticolo}.

vi) Suppose that $R_{G,\theta}$ has type $B_{2}$. We can realize
$R_{G,\theta}$ as in the case of $B_{l}$ with $l>2$. First suppose
that $\rho(D(\widetilde{X}))=\{\alpha^{\vee}_{1} \}$ and write
$v=-a\omega^{\vee}_{1}-b\omega^{\vee}_{2}$. Let
$\{\lambda_{1},\lambda_{2}\}$ be the dual basis of
$\{\alpha^{\vee}_{1},v\}$; we have
$\lambda_{1}=-\alpha_{2}-\frac{b}{a+2b}(-\alpha_{1}-2\alpha_{2})$
and $\lambda_{2}=\frac{1}{a+2b}(-\alpha_{1}-2\alpha_{2})$. We have
$a=1$ because $\frac{1}{2}(2\lambda_{1}-\lambda_{2})=
\frac{1}{2}(-2\alpha_{2}-\frac{1+2b}{a+2b}(-\alpha_{1}-2\alpha_{2}))$
must be $\frac{1}{2}\alpha_{1}$. On the other hand  $1+2b $ divides
2 because $\lambda_{2}$ is a weight, so $b=0$. Hence
$v=-\omega^{\vee}_{1}$; in particular $\widetilde{X}$ is complete
and $\chi(S)$ is the restricted root lattice because
$-\omega_{2}^{\vee}=-2\omega_{1}^{\vee}+\alpha^{\vee}_{1}$. We have
proved that there is only one smooth, simple embedding with Picard
number one and $\rho(D(\widetilde{X}))=\{\alpha^{\vee}_{1} \}$;
moreover it is complete. Thus there is no a  smooth, complete, non
simple embedding with Picard number one,  because otherwise it would
strictly contain a simple, open, smooth  $G$-subvariety
$\widetilde{X}$ with Picard number one and
$\rho(D(\widetilde{X}))=\{\alpha^{\vee}_{1} \}$. Finally, we
consider the case in which $X$ is simple and
$\rho(D(X))=\{\alpha^{\vee}_{2} \}$, so
$C(X)=cone(-\omega^{\vee}_{2},\alpha_{2}^{\vee})$. Because of the
completeness of $X$, one can study this case exactly as the case v)
of $C_{n}$. Observe that if $G/H$ is Hermitian, then we have to
consider only the first possibility, because in the second case the
Picard number would be two.

vii) Suppose now that $R_{G,\theta}$ has type $BC_{l}$. We can
realize $R_{G,\theta}$ as the set   $\{\pm(e_{i}\pm e_{j})\,|\,
i\neq j\}\cup\{\pm2e_{i}\}\cup\{\pm e_{i}\}\subset \mathbb{R}^{l}$;
we can also   choose  $\{e_{1}-e_{2}$,...,$e_{l-1}-e_{l}, e_{l}\}$
as basis of $R_{G,\theta}$. The weight lattice is equal to the root
lattice; moreover, $BC_{l}$ coincides with its dual root system in
this realization  as subset of $\mathbb{R}^{l}$ (but
$\frac{2}{(\alpha,\alpha)}\alpha$ can be different from $\alpha$).
We proceed as in the cases of $B_{l}$ and $C_{l}$. The $G$-variety
$X$ must be simple; hence $C(X)$ is the cone generated by $-C^{+}$
and $\rho(D(X))$, namely
$(cone(\alpha^{\vee}_{1},...,\alpha^{\vee}_{l-1},
-\omega^{\vee}_{1})$.   Observe that we exclude the case where $l=2$
and $\rho(D(\widetilde{X}))=\{\alpha^{\vee}_{2}\}$, because
$R_{L,\theta}$ would have type $BC_{1}$. The set
$\{\alpha^{\vee}_{1},...,\alpha^{\vee}_{l-1},$
$-\omega^{\vee}_{1}\}$ is a basis of the restricted coroot lattice
because $\omega^{\vee}_{1}=e_{1}$; now, one can prove the smoothness
of $X$ exactly as in the case of $B_{l}$.

viii) Suppose that $R_{G,\theta}$  has type $D_{l}$. We can realize
$R_{G,\theta}$  as the set  $\{\pm(e_{i}\pm e_{j})\,|\, i\neq
j\}\subset \mathbb{R}^{l}$; moreover we can choose
$\{e_{1}-e_{2}$,...,$e_{l-1}-e_{l},e_{l-1}+e_{l}\}$ as basis. In
particular, we can we identify the restricted roots with the
corresponding coroots. We can suppose, up to an automorphism of the
Dynkin diagram, that
$\rho(D(\widetilde{X}))=\{\alpha_{1},...,\alpha_{l-1}\}$.  Suppose
first that $(\alpha_{i},v)=-\delta_{i,1}$ for each $i<l$, so $v=
 \sum_{i=1}^{l-2}(ia+i-1)\alpha_{i}
+\frac{(l-2)a+l-3}{2}\alpha_{l-1}+\frac{la+l-1}{2}\alpha_{l}$ for an
appropriate constant $a$. The restricted root $\alpha_{l}$ belongs
to $\mathbb{Z}v\oplus\bigoplus_{i=1}^{l-1}\mathbb{Z}\alpha_{i}$,
thus its coordinate with respect to $v$, namely $\frac{2}{la+l-1}$,
is an integer. On the other hand $v=-\omega_{1}-b\omega_{l}$; we
obtain $2a=-b-2$ comparing the coordinates of $v$ in the bases of
the simple restricted roots, respectively of the fundamental
spherical weights. The previous facts imply that $b=0$, so $v=
-\omega_{1}$. We set $f_{i}=-e_{l-i+1}$ and
$\gamma_{i}=\alpha_{l-i}$, so that $-\omega_{1}=f_{l}$,
$\gamma_{0}=-f_{1}-f_{2}$ and $\gamma_{i}=f_{i}-f_{i+1}$ for each
$i=1,...,l-1$. Defining $\{\lambda_{1},...,\lambda_{l}\}$ as the
dual basis  of $\{\gamma_{1},...,\gamma_{l-1},v\}$, we have
$\lambda_{j}=\sum_{i=1}^{j}f_{i}$ for each $j$, thus
$\frac{1}{l}(l\lambda_{j}-j\lambda_{l})$ is the $j$-th fundamental
weight of $R_{L,\theta}$. Moreover the lattice $\chi_{*}(S)$ is
freely generated by
$\omega_{1},...,\omega_{l-2},\omega_{l-1}+\omega_{l},2\omega_{l}$.

Suppose now that  $(\alpha_{i},v)=-\delta_{i,l-1}$ for each $i<l$,
thus $v=-\omega_{l-1}-b\omega_{l}$. Moreover  $v$ has to be equal to
$a\sum_{i=1}^{l-2}i
\alpha_{i}+\frac{(l-2)a-1}{2}\alpha_{l-1}+\frac{la+1}{2}\alpha_{l}$
for an appropriate constant $a$.  We compare the coordinates of $v$
in the bases of the simple restricted roots, respectively of the
fundamental spherical weights, obtaining  $2a=-b-1$. On the other
hand, the restricted root $\alpha_{l}$ must belong to the lattice
generated by $\{\alpha_{1},...,\alpha_{l-1 },v\}$, so its coordinate
with respect to $v$, namely $\frac{2}{la+1}$, must be an integer;
hence $a=-\frac{1}{2}$,  $v=-\omega_{l-1}$ and   $l$ is equal either
to four or to six. If $l=4$ we can reduce this case to the one where
$C(\widetilde{X})=cone(\alpha_{1},\alpha_{2},\alpha_{3},-\omega_{1})$
by an automorphism of the Dynkin diagram. Finally consider the case
where $l=6$ and
$C(\widetilde{X})=cone(\alpha_{1},...,\alpha_{5},-\omega_{5})$. The
variety $\widetilde{X}$ must be complete because the simple
symmetric variety corresponding to
$(cone(\alpha_{1},...,\alpha_{4},\alpha_{6},-\omega_{5}),
\{D_{\alpha_{1}^{\vee}},...,D_{\alpha_{4}^{\vee}},D_{\alpha_{6}^{\vee}}\})$
is not smooth. But $-\omega_{3}$ is equal to $
-\frac{3}{2}\omega_{5} -\frac{1}{4}\alpha_{1} -\frac{1}{2}\alpha_{2}
-\frac{3}{4}\alpha_{3} +\frac{3}{4}\alpha_{5}$, so it does not
belong to $cone(\alpha_{1},...,\alpha_{5},-\omega_{5})$, a
contradiction.

We have proved that, if $l$ is different from four, then there is at
most one complete symmetric variety $X$ with the requested
properties and it is such that: i)
$\chi_{*}(S)=\bigoplus_{i=1}^{l-2}
\mathbb{Z}\omega_{i}\oplus\mathbb{Z}(\omega_{l-1}+\omega_{l}) \oplus
\mathbb{Z}2\omega_{l}$; ii) the corresponding colored fan is formed
by $(cone(\alpha_{1},...,\alpha_{l-1},-\omega_{1}),$
$\{D_{\alpha_{1}^{\vee}},...,D_{\alpha_{l-1}^{\vee}}\})$, by
$(cone(\alpha_{1},...,\alpha_{l-2},\alpha_{l},$ $-\omega_{1}),$
$\{D_{\alpha_{1}^{\vee}},...,$
$D_{\alpha_{l-2}^{\vee}},D_{\alpha_{l}^{\vee}}\})$ and by their
colored faces. We have to show that these combinatorial data define
a variety and that this variety is complete.   To verify that these
colored cones define a colored fan it is sufficient to prove that
the intersection of $cone(\alpha_{1},...,\alpha_{l-1},-\omega_{1})$
with $cone(\alpha_{1},...,\alpha_{l-2},\alpha_{l},-\omega_{1})$ is
$cone(\alpha_{1},...,\alpha_{l-2},-\omega_{1})$. Let
$\sum_{i=1}^{l-2}a_{i}\alpha_{i}+b\alpha_{l}-c\omega_{1}=
\sum_{i=1}^{l-2}(a_{i}-2b)\alpha_{i}-b\alpha_{l-1}+(-c+2b)\omega_{1}$
be a vector in the intersection, then $b$ is equal to zero. The
variety $X$ is complete because
$cone(\alpha_{1},...,\alpha_{l-1},-\omega_{1})$ $\cap\, -C^{+}=
cone(-\omega_{1},...,-\omega_{l-2},-\omega_{l-1}-\omega_{l},-\omega_{l})$
and $cone(\alpha_{1},...,\alpha_{l-2},\alpha_{l},-\omega_{1})\cap
-C^{+}=
cone(-\omega_{1},...,-\omega_{l-2},-\omega_{l-1}-\omega_{l},-\omega_{l-1})$.
If $l$ is equal to four, we can proceed in an analogous way.

If $R_{G,\theta}$  has type $E_{l}$, we number the simple roots so
that $\alpha_{1},\widehat{\alpha}_{2},\alpha_{3},...,\alpha_{l}$
generated a root system of type $A_{l-1}$; moreover we choose the
inner product so that $\alpha_{i}=\alpha_{i}^{\vee}$ for each $i$.
Thus $\rho(D(\widetilde{X}))$ has to be
$\{\alpha_{1},\widehat{\alpha}_{2},\alpha_{3},...,\alpha_{l}\}$ and
$\widetilde{X}$ must be  complete.

ix) Suppose that $R_{G,\theta}$  has type $E_{6}$.  Up to an
automorphism of the Dynkin diagram, we can suppose
$(\alpha_{i},v)=-\delta_{i,6}$ for each $i$ different from two. Thus
$v=\sum a_{i}\alpha_{i}=-a\omega_{2}-\omega_{6}$, where $a$ is a
positive integer. We compare the coordinates of $v$  in the bases of
the simple restricted roots, respectively of the fundamental
spherical weights, obtaining $a_{2}=-2a-1$. On the other hand,
$\alpha_{2}$ belongs to the lattice generated by $
\alpha_{1},\alpha_{3},...,\alpha_{6},v $, thus $(a_{2})^{-1}$ is an
integer; hence $v$ has to be $-\omega_{6}$. The weight
$-3\omega_{1}$ is equal to $-3\omega_{6}-2\alpha_{1}-\alpha_{3}
+\alpha_{5}+2\alpha_{6}$, so it does not belong to
$C(\widetilde{X})$. Thus $\widetilde{X}$ is not complete, a
contradiction.

x) Suppose that $R_{G,\theta}$  has type $E_{7}$ and write $v=\sum
a_{i}\alpha_{i}$. The coefficients $a_{i}$ belong to
$\frac{1}{2}\mathbb{Z}$ because $v$ is a weight. Moreover
$(a_{2})^{-1}$ is an integer, because $\alpha_{2}$ belongs to the
lattice generated by $\alpha_{1},\alpha_{3},...,\alpha_{7},v$, thus
$a_{2}$ is equal either to $\pm1$ or to $\pm\frac{1}{2}$. First
suppose that $(\alpha_{i},v)=-\delta_{i,1}$ for each $i$ different
from two, so $v=-\omega_{1}-a\omega_{2} $ for an appropriate
positive integer $a$. We compare the coordinates of $v$  in the
bases of the simple restricted roots, respectively of the
fundamental spherical weights, obtaining $a=\frac{2}{7}(-a_{2}-2)$;
thus $a$ is not integer, a contradiction. Finally,  suppose that
$(\alpha_{7},v)=-1$. We have only to study the  basis
$\{\lambda_{1},...,\lambda_{7}\}$, where
$\lambda_{1}=\omega_{1}-\frac{a_{1}}{a_{2}}\omega_{2}$,
$\lambda_{i}=\omega_{i+1}-\frac{a_{i+1}}{a_{2}}\omega_{2}$   for
$1<i<7$ and $\lambda_{7}=\frac{1}{a_{2}}\omega_{2}$. Comparing the
second coordinate of $\frac{1}{7}(7\lambda_{1}-\lambda_{7})$ (with
respect to the basis $\{\alpha_{1}, ...,\alpha_{7}\}$) with the one
of the first fundamental weight of $R_{L,\theta}$, we obtain
$7a_{1}=4a_{2}-1$. Thus $2(4a_{2}-1)$ is integral multiple of $7$, a
contradiction.

xi) Suppose that $R_{G,\theta}$  has type $E_{8}$ and write $v=\sum
a_{i}\alpha_{i}$. Notice that the root lattice of $R_{G,\theta}$ is
equal to the weight lattice. The integer $a_{2}$ is $\pm1$ because
$\{\alpha_{1},\alpha_{3},...,\alpha_{8},v\}$ is a basis of
$\chi_{*}(S)$. First, suppose that $v=-\omega_{1}-a\omega_{2} $ for
an appropriate positive integer $a$. We compare the coordinates of
$v$ in the two bases, respectively of simple restricted roots and
fundamental spherical weights, obtaining $a=\frac{1}{8}(-a_{2}-5)$;
so $a $ is not integer, a contradiction. Finally, suppose that
$(\alpha_{8},v)=-1$. We have only to study the indexed basis
$\{\lambda_{1},...,\lambda_{8}\}$ with
$\lambda_{1}=\omega_{1}-\frac{a_{1}}{a_{2}}\omega_{2}$,
$\lambda_{i}=\omega_{i+1}-\frac{a_{i+1}}{a_{2}}\omega_{2}$ for
$1<i<8$ and $\lambda_{8}=\frac{1}{a_{2}}\omega_{2}$. Comparing the
second coordinate of $\frac{1}{8}(8\lambda_{1}-\lambda_{8})$ (with
respect to the basis $\{\alpha_{1},...,\alpha_{8}\}$)  with the one
of the first fundamental weight of $R_{L,\theta}$, we obtain
$8a_{1}=5a_{2}-1$. On the other hand, $5a_{2}-1$ cannot be a
integral multiple of $8$, a contradiction.

xii) If $R_{G,\theta}$  has type $F_{4}$, then $R_{L,\theta}$ cannot
have type $A_{3}$ and there are no simple smooth varieties with
Picard number  one.

xiii) Suppose that $R_{G,\theta}$  has type $G_{2}$ (and assume
$\alpha_{1}$ short).  Write
$v=-a\omega^{\vee}_{1}-b\omega^{\vee}_{2}$ and observe that the
weight lattice coincides with the root lattice. First, suppose that
$D(X)$ contains $D_{\alpha_{1}^{\vee}}$, hence $X$ contains a simple
smooth subvariety $\widetilde{X}$ with
$C(\widetilde{X})=cone(\alpha_{1}^{\vee},v)$. Let
$\{\lambda_{1},\lambda_{2}\}$ be the dual basis of
$\{\alpha_{1}^{\vee},v\}$. We have
$\lambda_{1}=-\frac{1}{3}\alpha_{2}+\frac{b}{9a+6b}(3\alpha_{1}+2\alpha_{2})$
and $\lambda_{2}= -\frac{1}{3a+2b}(3\alpha_{1}+2\alpha_{2})$. Thus
$3a+2b$ must  divide 2 and 3, so it must be 1, a contradiction. Thus
$X$ is simple and $C(X)=cone(\alpha_{2}^{\vee},-\omega^{\vee}_{2})$.
The dual basis of $\{\alpha_{2}^{\vee},-\omega_{2}^{\vee}\}$ is $\{
-\alpha_{1},-\alpha_{2}-2\alpha_{1}\}$. Moreover
$\frac{1}{2}(2\lambda_{1}-\lambda_{2})= \frac{1}{2} \alpha_{2}$ is
the fundamental weight of $R_{L,\theta}$. Thus the variety
associated to
$(cone(\alpha_{2}^{\vee},-\omega^{\vee}_{2}),D_{\alpha_{2}^{\vee}})$
is smooth.

\begin{lem}\label{proj}
Every smooth completion $X$ of $G/H$ with Picard number one  is
projective.
\end{lem}

{\em Proof.} It is sufficient to consider the varieties which are
not simple. There are exactly two maximal colored cones, say
$(\varsigma_{1},I_{1})$ and  $(\varsigma_{2},I_{2})$; moreover there
are exactly two colors, say $D_{1}\in I_{1}$ and $D_{2}\in I_{2}$,
which do not belong to $I_{1}\cap I_{2}$. We claim that
$D_{1}+D_{2}$ is an ample divisor. Indeed, let $\varphi$ be the
function over $\varsigma_{1}\cup \varsigma_{2}$ corresponding to
$D_{1}+D_{2}$ (see Theorem 3.1 in \cite{Br1}) and let $l_{i}$ be the
linear function which coincides with $\varphi$ over $\varsigma_{i}$.
The cone $\varsigma_{i}$ is generated by $\rho(D_{i})$ and $
\varsigma_{1}\cap \varsigma_{2}$; moreover $\rho(D_{1})$ and
$\rho(D_{2})$ are permutated by the reflection with respect to the
hyperplane generated by  $ \varsigma_{1}\cap \varsigma_{2}$. Thus
$l_{1}(\rho(D_{2}))=-1<1=\varphi(\rho(D_{2}))$ and
$l_{2}(\rho(D_{1}))=-1<1=\varphi(\rho(D_{1}))$, because
$span_{\mathbb{R}}(\varsigma_{1}\cap \varsigma_{2})$ is the kernel
of both $l_{1}$ and $l_{2}$. Hence $\varphi$ is strictly convex over
the colored fan of $X$. Moreover $D(X)=D(G/H)$, therefore
$D_{1}+D_{2}$ is ample by Theorem 3.3 in \cite{Br1}. $\square$

We have proved the following theorem:

\begin{thm}\label{classif} Let $G$ be a semisimple simply connected group and let
$G/H$ be a homogeneous symmetric variety.  Suppose that there is a
smooth, complete embedding $X$ of $G/H$ with Picard number one.
Then:
\begin{itemize}
\item Given $G/H$, there is, up to  equivariant isomorphism, at most one  embedding with
the previous properties.

\item The symmetric variety $X$  is projective.

\item The number of colors of $G/H$ is equal to the rank $l$ of $G/H$; in
particular there are no exceptional roots. If, in addition, $l$ is
at least three, then  $G/H$ is not Hermitian.
\item The restricted root system $R_{G,\theta}$ is irreducible or has type $A_{1}\times
A_{1}$.
\item We have two possibilities:
\begin{enumerate}
\item[(i)]  $X$ is simple and $D(X)$ has cardinality $l-1$.
\item[(ii)]  $X$ contains two closed orbits, $D(X)$ is equal to
$D(G/H)$ and $H$ has index two or three in $N_{G}(G^{\theta})$.
\end{enumerate}
In particular, $X$ is simple if $H=N_{G}(G^{\theta})$.
\item We have the following classification
depending on the type of the restricted root system $R_{G,\theta}$ :
\begin{enumerate}
\item[(i)] If $R_{G,\theta}$  has type $A_{1}\times A_{1}$, then $\chi(S)$ has basis
$\{2\omega_{1},\omega_{1}+\omega_{2}\}$; in particular $H$ has index
two in $N_{G}(G^{\theta})$. Moreover, $X$ has two closed orbits; the
maximal colored cones of the colored fan of $X$ are
$(cone(\alpha^{\vee}_{1},-\omega^{\vee}_{1}-\omega^{\vee}_{2}),$ $
\{D_{\alpha^{\vee}_{1}}\})$ and
$(cone(\alpha^{\vee}_{2},-\omega^{\vee}_{1}-\omega^{\vee}_{2}),$ $
\{D_{\alpha^{\vee}_{2}}\})$.

\item[(ii)] If $l=1$, then $G/H$ can   be isomorphic neither to
$SL_{n+1}/S(L_{1}\times L_{n})$, nor to $SL_{2}/SO_{2}$. With such
hypothesis, $G/H$ has a  unique non trivial embedding which is
simple, projective, smooth and with Picard number 1.
\item[(iii)] If $R_{G,\theta}$  has type $A_{l}$ with $l>1$, we have the following
possibilities:
\begin{itemize}
\item[(a)] $H=N_{G}(G^{\theta})$ and $X$ is simple. In this case $X$ is
associated either to the colored cone
$(cone(\alpha^{\vee}_{1},...,\alpha^{\vee}_{l-1},-\omega^{\vee}_{1}),
\{D_{\alpha^{\vee}_{1}},...,D_{\alpha^{\vee}_{l-1}}\})$ or to the
one
$(cone(\alpha^{\vee}_{2},...,\alpha^{\vee}_{l},-\omega^{\vee}_{l}),
\{D_{\alpha^{\vee}_{2}},...,D_{\alpha^{\vee}_{l }}\})$;
\item[(b)] $H=G^{\theta}$ and $l=2$. In this case $X$ has two closed orbits. The
maximal colored cones of the colored fan of $X$ are
$(cone(\alpha^{\vee}_{1},-\omega^{\vee}_{1}-\omega^{\vee}_{2}),
\{D_{\alpha^{\vee}_{1}} \})$ and
$(cone(\alpha^{\vee}_{2},-\omega^{\vee}_{1}-\omega^{\vee}_{2}),
\{D_{\alpha^{\vee}_{2}}\})$.
\end{itemize}
\item[(iv)] If $R_{G,\theta}$  has type $B_{2}$, then $X$ is simple and  we have the following
possibilities:

\begin{itemize}
\item[(a)] $H=N_{G}(G^{\theta})$ and  $X$ is
associated   to the colored cone
$(cone(\alpha^{\vee}_{1},-\omega^{\vee}_{1}),$ $
\{D_{\alpha^{\vee}_{1}}\})$;
\item[(b)] $H=G^{\theta}$ and $X$ is
associated   to the colored cone
$(cone(\alpha^{\vee}_{2},-\omega^{\vee}_{2}),$
$\{D_{\alpha^{\vee}_{2}}\})$. Moreover  $G/H$ cannot be Hermitian.
\end{itemize}

\item[(v)] If $R_{G,\theta}$  has type $B_{l}$ with $l>2$, then $H=N_{G}(G^{\theta})$,  $X$ is
simple and  is associated  to the colored cone
$(cone(\alpha^{\vee}_{1},...,\alpha^{\vee}_{l-1},-\omega^{\vee}_{1}),
\{D_{\alpha^{\vee}_{1}},...,$ $D_{\alpha^{\vee}_{l-1}}\})$.

\item[(vi)] If $R_{G,\theta}$  has type $C_{l}$,   then $H= G^{\theta} $, $X$ is simple and
corresponds to the colored cone
$(cone(\alpha^{\vee}_{1},...,\alpha^{\vee}_{l-1},-\omega^{\vee}_{1}),
\{D_{\alpha^{\vee}_{1}},...,D_{\alpha^{\vee}_{l-1}}\})$. Moreover
$G/H$ cannot be Hermitian.

\item[(vii)] If $R_{G,\theta}$  has type $BC_{l}$ with $l>1$, then $H=N_{G}(G^{\theta})= G^{\theta}$,  $X$ is
simple and corresponds  to the colored cone
$(cone(\alpha^{\vee}_{1},...,\alpha^{\vee}_{l-1},-\omega^{\vee}_{1}),$
$ \{D_{\alpha^{\vee}_{1}},...,$ $D_{\alpha^{\vee}_{l-1}}\})$.
\item[(viii)] If $R_{G,\theta}$  has type $D_{l}$ with $l>4$, then $\chi_{*}(S)$ is freely generated by
$\omega^{\vee}_{1},...,\omega^{\vee}_{l-2},\omega^{\vee}_{l-1}+\omega^{\vee}_{l},2\omega^{\vee}_{l}$;
in particular $H$ has index two in $N_{G}(G^{\theta})$. $X$ has two
closed orbits; the maximal colored cones of the colored fan of $X$
are
$(cone(\alpha^{\vee}_{1},...,\alpha^{\vee}_{l-1},-\omega^{\vee}_{1}),$
$ \{D_{\alpha^{\vee}_{1}},...,D_{\alpha^{\vee}_{l-1}}\})$ and
$(cone(\alpha^{\vee}_{1},..., \alpha^{\vee}_{l-2},$
$\alpha^{\vee}_{l}, -\omega^{\vee}_{1}),$
$\{D_{\alpha^{\vee}_{1}},...,$ $D_{\alpha^{\vee}_{l-2}},$
$D_{\alpha^{\vee}_{l }}\})$.
\item[(ix)] If $R_{G,\theta}$  has type $D_{4}$, then $H$ has index two in $N_{G}(G^{\theta})$ and $X$
has two closed orbits. If
$\chi_{*}(S)=\mathbb{Z}\omega^{\vee}_{i}\oplus
\mathbb{Z}\omega^{\vee}_{2}\oplus
\mathbb{Z}(\omega^{\vee}_{j}+\omega^{\vee}_{k})\oplus
\mathbb{Z}2\omega^{\vee}_{k}$, then the maximal colored cones of the
colored fan of  $X$ are
$(cone(\alpha^{\vee}_{i},\alpha^{\vee}_{2},\alpha^{\vee}_{j},-\omega^{\vee}_{i}),$
$
\{D_{\alpha^{\vee}_{i}},D_{\alpha^{\vee}_{2}},D_{\alpha^{\vee}_{j}}\})$
and
$(cone(\alpha^{\vee}_{i},\alpha^{\vee}_{2},\alpha^{\vee}_{k},-\omega^{\vee}_{i}),$
$ \{D_{\alpha^{\vee}_{i}},D_{\alpha^{\vee}_{2}},$
$D_{\alpha^{\vee}_{k}}\})$.
\item[(x)]
The type of $R_{G,\theta}$  cannot be $E_{6}$, $E_{7}$, $E_{8}$ or
$F_{4}$.

\item[(xi)] If $R_{G,\theta}$  has type $G_{2}$ then $H=N_{G}(G^{\theta})= G^{\theta}$,  $X$ is
simple and  is associated  to the colored cone
$(cone(\alpha^{\vee}_{2},-\omega^{\vee}_{2}),
\{D_{\alpha^{\vee}_{2}}\})$.

\end{enumerate}
\end{itemize}
\end{thm}

\section*{Acknowledgments}

We would like to thank M. Brion  for many very helpful suggestions.
I would also thank the referee because his accurate reports are been
very useful to me.


\end{document}